\documentclass[a4paper]{style_used}

\usepackage{tikz}
\usepackage{subcaption}
\usepackage{amssymb}
\usepackage{graphicx}
\usepackage{amsmath}
\usepackage{amsthm}
\usepackage{accents}
\usepackage{mathrsfs}

\newcommand{\R}{\mathbb{R}}

\newcommand{\f}[1]{\mathbf{#1}}
\newcommand{\cupdot}{\,\dot{\cup}\,}
\newcommand{\Du}{\partial_{1}}
\newcommand{\Dv}{\partial_{2}}
\newcommand{\ev}{ \accentset{\star}}

\DeclareMathOperator{\Span}{span}

\title[Isogeometric analysis with $C^1$ functions]{Isogeometric analysis with $C^1$ functions \\ on unstructured quadrilateral meshes}

\author[M. Kapl]{\firstname{Mario} \lastname{Kapl}}
\address{Johann Radon Institute for Computational and Applied Mathematics, Austrian Academy of Sciences, Austria}
\email{mario.kapl@ricam.oeaw.ac.at}

\author[G. Sangalli]{\firstname{Giancarlo} \lastname{Sangalli}}
\address{Dipartimento di Matematica ``F. Casorati'', Universit\`a degli Studi di Pavia, and IMATI-CNR, Pavia, Italy}
\email{giancarlo.sangalli@unipv.it}

\author[T. Takacs]{\firstname{Thomas} \lastname{Takacs}}
\address{Institute of Applied Geometry, Johannes Kepler University Linz, Austria}
\email{thomas.takacs@jku.at}



\begin{document}

\begin{abstract}
In the context of isogeometric analysis, globally $C^1$ isogeometric  spaces over unstructured quadrilateral meshes allow the direct solution of fourth order partial differential 
equations on complex geometries via their Galerkin discretization. The design of such smooth spaces has been intensively studied in the last five years, in particular for the case of 
planar domains, and is still task of current research. In this paper, we first give a short survey of the developed methods and especially focus on the approach~\cite{KaSaTa17c}. There, 
the construction of a specific $C^1$ isogeometric spline space for the class of so-called analysis-suitable $G^1$ multi-patch parametrizations  is presented. This particular 
class of parameterizations comprises exactly those multi-patch geometries, which ensure the design of $C^1$ spaces with optimal approximation properties, and allows the representation 
of complex planar multi-patch domains. We present known results in a coherent framework, and also extend the construction to parametrizations that are not analysis-suitable $G^1$ by 
allowing  higher-degree splines in the neighborhood of the extraordinary vertices and edges. Finally, we present numerical tests that illustrate the behavior of the proposed method 
on representative examples.
\end{abstract}

\maketitle


\section{Introduction}

Isogeometric analysis (IgA) was introduced in~\cite{HuCoBa04} as a framework for numerically solving partial differential equations (PDEs). The basic idea is to bridge the 
gap between geometric modeling (that is, Computer-Aided Design) and numerical analysis (that is, Finite Element Analysis) by using the same (rational) spline function space for 
representing the geometry of the computational domain and for describing the solution of the PDE (cf.~\cite{ANU:9260759,CottrellBook,HuCoBa04}). In contrast to finite 
elements, IgA allows a simple integration of smooth discretization spaces for the numerical simulation. While the design of such smooth spaces is trivial for 
single patch geometries, it is a challenging task for the case of multi-patch or manifold geometries. 

The scope of this paper is to give a survey of the different existing methods for the construction of strongly enforced $C^1$ isogeometric spline spaces over planar, 
unstructured quadrilateral meshes (that is, planar multi-patch geometries composed of quadrilateral patches with possibly extraordinary vertices), see Section~\ref{sec:design}, 
with a special focus on the approach~\cite{KaSaTa17c}, see Section~\ref{sec:C1ASG1space}. The common goal of the developed techniques is to generate isogeometric spline spaces, 
which are not only exactly $C^1$-smooth within the single patches but also across the patch interfaces. The design of the $C^1$ isogeometric spaces is mainly based on the observation 
that an isogeometric function is $C^1$-smooth if and only if the associated graph surface is $G^1$-smooth (that is, geometric continuous of order~$1$), cf.~\cite{Pe15}. 
The resulting global $C^1$-smoothness of the spaces then enables the solution of fourth order PDEs just via its weak form using a standard Galerkin discretization, see 
for example  \cite{BaDe15,CoSaTa16,KaBuBeJu17,KaViJuBi15,NgKaPe15,TaDe14} 
for the biharmonic equation,  \cite{ABBLRS-stream,benson2011large,kiendl-bazilevs-hsu-wuechner-bletzinger-10,kiendl-bletzinger-linhard-09,KiHsWuRe15} for the 
Kirchhoff-Love shell formulation, \cite{gomez2008isogeometric,GoCaHu09,LiDeEvBoHu13} for the Cahn-Hilliard equation and \cite{gradientElast2011,KhakaloNiiranenC1} for plane 
problems of first strain gradient elasticity.

A further possible strategy to impose $C^1$-smoothness across the interfaces of general multi-patch geometries is the use of 
subdivision surfaces (e.g., ~\cite{BHU10b,CiOrSch00,CiScAnOrSch02,JuMaPeRu16,Peters2,RiAuFe16,ZhSaCi18}). The surfaces are recursively generated via refinement 
schemes, and are described in the limit as the collection of infinitely many polynomial patches, e.g., in the case of Catmull-Clark subdivision 
of bicubic patches. We refer to~\cite{Subdivision_surfaces} for further reading. Challenges of dealing with subdivision surfaces in IgA include the need for special techniques for the 
numerical integration~\cite{JuMaPeRu16} and the often reduced approximation power in the neighborhood of extraordinary vertices~\cite{Peters2}.

Instead of enforcing the $C^1$-continuity across the patch interfaces in a strong sense, the $C^1$-smoothness could also be achieved by coupling the neighboring patches in a weak sense. 
We do not cover here this approach, which is typically based on adding penalty terms to the weak formulation of the PDE~\cite{ApBrWuBl15,Guo2015881},  or using Lagrange 
multipliers~\cite{ApBrWuBl15,BenvenutiPhD}. These techniques are applicable to quite general multi-patch geometries with even non-matching meshes,  but at the cost of obtaining an 
approximate $C^1$ solution.  Moreover, the formulation of the problem, and as a result the system matrix, have to be adapted accordingly.

The outline of this paper is as follows. Section~\ref{sec:design} describes the state of the art of such constructions over planar multi-patch domains. 
Section~\ref{sec:multi-geometries} further discusses the case of  multi-patch parameterizations which are regular (that is, non-singular) and $C^0$ at the patch interfaces, 
and specifically describes so-called  analysis-suitable $G^1$ parametrizations. In this setting  we review the construction of $C^1$ isogeometric spaces, and their properties, 
in Section \ref{sec:C1ASG1space}. Section \ref{sec:generalization} extends the construction beyond  analysis-suitable $G^1$ parametrizations, by allowing higher-order 
splines around the extraordinary vertices. Numerical evidence of the optimal order of convergence of the proposed method is reported in Section \ref{sec:numerical_example}.

\section{The design of $C^1$ isogeometric spaces} \label{sec:design}

We give an overview of existing strategies for the design of strongly enforced $C^1$ isogeometric spline spaces over unstructured quadrilateral meshes on planar domains. 
Such quadrilateral meshes can be understood in the context of multi-patch domains or spline manifolds. In the language of manifolds, we have given local charts, which usually overlap. 
The global smoothness is then determined by the smoothness within every chart. In the multi-patch framework, the patches do not overlap but share common interfaces. Hence, the smoothness 
is determined by the smoothness within the patches as well as the smoothness across interfaces. Since most CAD systems are built upon multi-patch structures we focus on this point of 
view.

The different techniques can be roughly classified into three approaches depending on the smoothness of the underlying parameterization~$\f{F}$ for the multi-patch domain~$\Omega$.

\paragraph{Multi-patch parameterizations which are $C^1$-smooth everywhere:}
In this setting, the parameterization of the multi-patch domain is assumed to be $C^1$-smooth everywhere. This immediately leads to a singularity appearing at every extraordinary 
vertex (EV). Consequently, the isogeometric functions are then $C^1$ everywhere away from the EV and possibly only $C^0$ at the EV, due to the singularity. To circumvent this issue, 
one has to enforce additional $G^1$ constraints at the EVs. This technique is based on the use of specific degenerate patches (e.g. D-patches~\cite{Re97}) in the neighborhood of an 
extraordinary vertex. These patches are obtained by collapsing some control points into one point and by a special configuration of some of the remaining control points which guarantee 
that the surface is $G^1$-smooth at the EV despite having a singularity there. The same approach allows to construct isogeometric functions that are  $C^1$ everywhere 
on the multi-patch domain. Examples of this method are~\cite{NgPe16,ToSpHu17b,ToSpHu17}. While the constructions~\cite{NgPe16,ToSpHu17} are restricted to bicubic splines, 
the methodology~\cite{ToSpHu17b} can be applied to bivariate splines of arbitrary bidegree~$(p,p)$ with $p \geq 3$. All three techniques can be used to construct sequences of nested 
isogeometric splines spaces.

A similar approach is to use subdivision surfaces to represent the isogeometric spaces. In subdivisions, the surface around an EV  is composed of an infinite sequence of spline rings, 
where every ring is (at least) $C^1$ smooth. The shape of the surface at the EVs is guided by the refinement rules of the control mesh. For example, for Catmull-Clark subdivision this 
procedure generates a surface that is $C^2$ smooth everywhere and $G^1$ at the EVs. 
However, the approach suffers from a lack of approximation power near the EV. See~\cite{BHU10b,Peters2,RiAuFe16,ZhSaCi18}, where subdivision based isogeometric analysis was studied.
 
\paragraph{Multi-patch parameterizations which are $C^1$-smooth except in the vicinity of an extraordinary vertex:}
The core idea is to construct a parameterization of the multi-patch domain which is $C^1$-smooth in the regular regions of the mesh and only $C^0$-smooth in a neighborhood of the 
extraordinary vertices (see e.g. \cite{BuJuMa16}). To obtain a globally $C^1$ isogeometric space, a $G^1$ surface construction is employed in the neighborhood of the EV. The same 
construction is used to generate the $C^1$ isogeometric functions over the multi-patch domain. This construction leads in general to a multi-patch surface which is even $C^{p-1}$-smooth 
away from an extraordinary vertex. One possibility is to use so-called G-splines, see \cite{Re95}. To obtain surfaces of good shape also in the vicinity of an extraordinary vertex, the 
$G^1$-smoothness is obtained by using a suitable surface cap, which requires a slightly higher degree than the surrounding $C^1$ spline surface. The resulting smooth surfaces are then 
used to construct $C^1$ isogeometric splines spaces, but which are in general not nested, see e.g.~\cite{Pe15-2,NgKaPe15}. The method~\cite{NgKaPe15} employs the surface 
construction~\cite{KaPe15A}, which is based on a biquadratic $C^1$ spline surface and on a bicubic or biquartic $G^1$ surface cap depending on the valency of the corresponding 
extraordinary vertex. The paper~\cite{Pe15-2} presents a new surface construction, where bicubic splines are complemented by biquartic splines in the neighborhood of extraordinary 
vertices. The methodology~\cite{KaPe18} can be seen as an extension of the above two techniques, and allows the construction of nested $C^1$ isogeometric spline spaces for a finite 
number of refinement steps. 

\paragraph{Multi-patch parameterizations which are $C^0$-smooth at all interfaces:}

The main idea is to consider multi-patch parametrizations that are everywhere regular (non-singular) but only $C^0$ at the patch interfaces, and then construct $C^1$ isogeometric 
spaces over them. Again, the key issue for application to isogeometric analysis is to guarantee good approximation properties of these spaces. 
In \cite{BlMoVi17,mourrain2015geometrically} the authors gave dimension formulas for meshes of arbitrary topology, consisting of quadrilateral polynomial patches and specific 
macro-elements. The techniques in~\cite{BlMoVi17,mourrain2015geometrically} work for splines of general bidegree~$(p,p)$ for some large enough $p$, and generate the $C^1$ basis 
functions by analyzing the module of syzygies of specific polynomial/spline functions. Extending from polynomials to general spline patches, dimensions were given and basis functions 
were constructed for bilinear two-patch domains in \cite{KaViJuBi15}. In \cite{CoSaTa16} the reproduction properties of the $C^1$-smooth subspaces along an interface were studied 
for arbitrary B-spline patches. From the presented results, bounds for the dimension of the $C^1$-smooth subspaces of arbitrary geometries can be derived. Moreover, the specific 
class of analysis-suitable $G^1$ parametrizations was identified. In the last few years, a number of methods were developed which follow this approach, and which allow in most 
cases the design of nested $C^1$ isogeometric spline spaces. These techniques use particular classes of $C^0$ regular multi-patch parameterizations to  obtain $C^1$ isogeometric 
spaces with good/optimal approximation properties: 

\begin{itemize}
 \item \emph{(Mapped) bilinear multi-patch parameterizations (e.g.~\cite{BeMa14,KaBuBeJu17,KaViJuBi15,Matskewich-PhD}):} The aim is to explore $C^1$ isogeometric spaces 
 over bilinear or mapped bilinear multi-patch geometries. The methods~\cite{BeMa14,Matskewich-PhD} study the spaces of biquintic and for some specific cases also 
 biquartic $C^1$ 
 isogeometric B\'{e}zier functions, and generate basis functions which are implicitly given by minimal determining sets (cf.~\cite{LaSch07}) for the involved B\'{e}zier 
 coefficients. In contrast to the other techniques using $C^0$ regular multi-patch parameterizations, the resulting spaces are not nested.
 In \cite{KaBuBeJu17,KaViJuBi15}, the case of bicubic and biquartic $C^1$ spline elements is considered. The resulting basis functions are explicitly given by simple 
 formulae and possess a small local support. While the paper~\cite{KaViJuBi15} deals with the case of two patches, the work in~\cite{KaBuBeJu17} is an extension to the 
 multi-patch case. 
 \item \emph{General analysis-suitable parameterizations (e.g. \cite{CoSaTa16,KaSaTa17a,KaSaTa17c}):} The previous strategy was based on simple geometries 
 such as bilinear parameterizations. The following approach uses a more general class of geometries, called analysis-suitable $G^1$ parameterizations, 
 cf.~\cite{CoSaTa16} and Section~\ref{subsec:analysis-suitable}, which includes the previous types of geometries. The class of analysis-suitable $G^1$ 
 parameterizations contains exactly those geometries which allow the design of $C^1$ isogeometric spline spaces with optimal approximation properties. 
 In \cite{KaSaTa17a}, the space of $C^1$ isogeometric spline functions over analysis-suitable $G^1$ two-patch geometries was analyzed. The developed method is 
 applicable to splines of general bidegree~$(p,p)$ with $p \geq 3$ and patch regularity $1 \leq r \leq p-2$, and constructs simple, explicitly given basis functions 
 with a small local support. The work~\cite{KaSaTa17c} extends the construction~\cite{KaSaTa17a} to the case of analysis-suitable $G^1$ multi-patch parameterizations 
 and will be discussed in detail in Section~\ref{sec:C1ASG1space}  
 
 \item \emph{Non-analysis-suitable parametrizations with elevated degree at the interfaces (e.g.~\cite{ChAnRa18}):} Instead of using a particular class of multi-patch parameterizations 
 for the multi-patch domain, the method~\cite{ChAnRa18} increases locally along the patch interfaces the bidegree of the $C^1$ isogeometric spline functions to get spaces with good 
 approximation properties. The constructed $C^1$ basis functions are implicitly given by means of minimal determining sets for the spline coefficients and possess in 
 general large supports over one or more entire patch interfaces. This approach is a direct consequence of the results presented in~\cite{CoSaTa16} and extends the ideas of Theorems 
 1 and 3 therein. We will present further theoretical foundations in this paper.
\end{itemize}

\paragraph{Polar configurations:} 
This approach is outside of the framework of unstructured multi-patch domains, as it can be interpreted as a regular mesh in polar coordinates. However, many ideas to study and 
construct smooth polar configurations can be carried over to extraordinary vertices. The approach is based on the use of polar splines to model the domains and to construct $C^1$ 
isogeometric spline spaces over these domains, see e.g.~\cite{Peters2,ToSpHiHu16}. The method~\cite{Peters2} employs the polar surface construction of~\cite{MyKaPe08}, which 
generates a bicubic $C^1$ polar spline surface. In~\cite{ToSpHiHu16}, a novel polar spline technology for splines of arbitrary bidegree~$(p,p)$ is presented. It can be used to 
generate polar spline surfaces which are $C^s$-smooth ($s \geq 0$) everywhere except at the polar point where the resulting surface is discontinuous. Furthermore, 
it was shown that this surface construction can be used to generate globally $C^k$ isogeometric function spaces.

\section{Multi-patch geometries and their analysis-suitable $G^1$ parameterization} \label{sec:multi-geometries}
We now focus on multi-patch parameterizations which are regular and $C^0$ at all interfaces. After some preliminaries and notation on  multi-patch geometries, we will consider one 
specific class of geometries, called analysis-suitable $G^1$-multi-patch parameterizations (cf.~\cite{CoSaTa16}), which will be used throughout the paper. Furthermore, we will use 
in the following a slightly adapted notation and definitions mainly based on our work developed in~\cite{KaSaTa17c} and \cite{KaSaTa17b}. 

\subsection{Multi-patch domain}

Let $p\geq 3$, $1 \leq r \leq p-2$ and $n \geq 1$. We denote by $\mathcal{S}^{p,r}_{h}$ the univariate spline space of degree~$p$ and continuity $C^{r}$ on the 
parameter domain~$[0,1]$ possessing the uniform open knot vector 
\begin{equation*}  
(\underbrace{0,\ldots,0}_{(p+1)-\mbox{\scriptsize times}},
\underbrace{\textstyle \frac{1}{n}\ldots ,\frac{1}{n}}_{(p-r) - \mbox{\scriptsize times}}, 
\underbrace{\textstyle \frac{2}{n},\ldots ,\frac{2}{n}}_{(p-r) - \mbox{\scriptsize times}},\ldots, 
\underbrace{\textstyle \frac{n-1}{n},\ldots ,\frac{n-1}{n}}_{(p-r) - \mbox{\scriptsize times}},
\underbrace{1,\ldots,1}_{(p+1)-\mbox{\scriptsize times}})
\end{equation*}
with the mesh size~$h=\frac{1}{n}$, and by $\mathcal{S}_{h}^{\f{p},\f{r}}$ with $\f{p}=(p,p)$ and $\f{r}=(r,r)$ the corresponding bivariate tensor-product spline 
space~$\mathcal{S}^{p,r}_{h} \otimes \mathcal{S}^{p,r}_{h}$ on the parameter domain~$[0,1]^{2}$. In addition, let $b_{j}$, $j=0,\ldots,N-1$, with 
$N =p+(n-1)(p-r)+1$, be the B-splines of 
$\mathcal{S}_{h}^{p,r}$, and let $b_{\f{j}}$, $\f{j}=(j_1,j_2) \in \{0,\ldots,N-1 \}^{2}$, be the tensor-product B-splines of $\mathcal{S}_{h}^{\f{p},\f{r}}$, that is,
\[
 b_{(j_1,j_2)}(\xi_1,\xi_2) = b_{j_1}(\xi_1)b_{j_2}(\xi_2).
\]

Consider an open domain~$\Omega \subset \R^2$, which is given as the union of quadrilateral patches~$\Omega^{(i)}$, $i \in \mathcal{I}_{\Omega}$, 
interfaces~$\Sigma^{(i)}$, $i \in \mathcal{I}^{\circ}_{\Sigma}$, and inner vertices~$\f{x}^{(i)}$, $i \in \mathcal{I}^{\circ}_{\chi}$, that is,
\[
 \Omega = \left( \bigcup_{i \in \mathcal{I}_{\Omega}}\Omega^{(i)}  \right) \cup \left( \bigcup_{i \in \mathcal{I}^{\circ}_{\Sigma}} \Sigma^{(i)}\right) 
 \cup \left( \bigcup_{i \in \mathcal{I}_{\chi}^{\circ}} \f{x}^{(i)} \right).
\]
We assume that all patches are mutually disjoint and that no hanging nodes exist. The boundary~$\Gamma$ of $\Omega$, that is, $\Gamma = \partial \Omega$, is given as the 
collection of boundary edges~$\Sigma^{(i)}$, $i \in \mathcal{I}^{\Gamma}_{\Sigma}$, and boundary vertices~$\f{x}^{(i)}$, $i \in \mathcal{I}_{\chi}^{\Gamma}$, that is,
\[
 \Gamma = \left( \bigcup_{i \in \mathcal{I}_{\Sigma}^{\Gamma}} \Sigma^{(i)} \right) \cup \left( \bigcup_{i \in \mathcal{I}^{\Gamma}_{\chi}} \f{x}^{(i)} \right) .
\]
In addition, we assume that $\mathcal{I}_{\Sigma}^{\circ} \cap \mathcal{I}_{\Sigma}^{\Gamma} = \emptyset$ and $\mathcal{I}_{\f{x}}^{\circ} \cap 
\mathcal{I}_{\f{x}}^{\Gamma} = \emptyset$, and denote by $\mathcal{I}_{\Sigma}$ and $\mathcal{I}_{\chi}$ the index sets 
$\mathcal{I}_{\Sigma} = \mathcal{I}_{\Sigma}^{\circ} \cupdot \mathcal{I}_{\Sigma}^{\Gamma}$ and 
$\mathcal{I}_{\chi} = \mathcal{I}_{\chi}^{\circ} \cupdot \mathcal{I}_{\chi}^{\Gamma}$, respectively. 

Each quadrilateral patch $\Omega^{(i)}$ is the open image of a bijective and regular geometry mapping
\[
\f{F}^{(i)}: [0,1]^2 \rightarrow
        \overline{\Omega^{(i)}}\subset \R^2 , 
\]
with $\f{F}^{(i)} \in \mathcal{S}_{h}^{\f{p},\f{r}} \times \mathcal{S}_{h}^{\f{p},\f{r}}$. We denote by $\f{F}$ the resulting multi-patch geometry of $\Omega$ consisting 
of the single geometry mappings~$\f{F}^{(i)}$, $i \in \mathcal{I}_{\Omega}$.

\subsection{Analysis-suitable $G^1$ parameterization: definition} \label{subsec:analysis-suitable}

\begin{figure}[ht!]
\centering
 \begin{picture}(250,125)
    \put(0,0){\includegraphics[width=.45\textwidth]{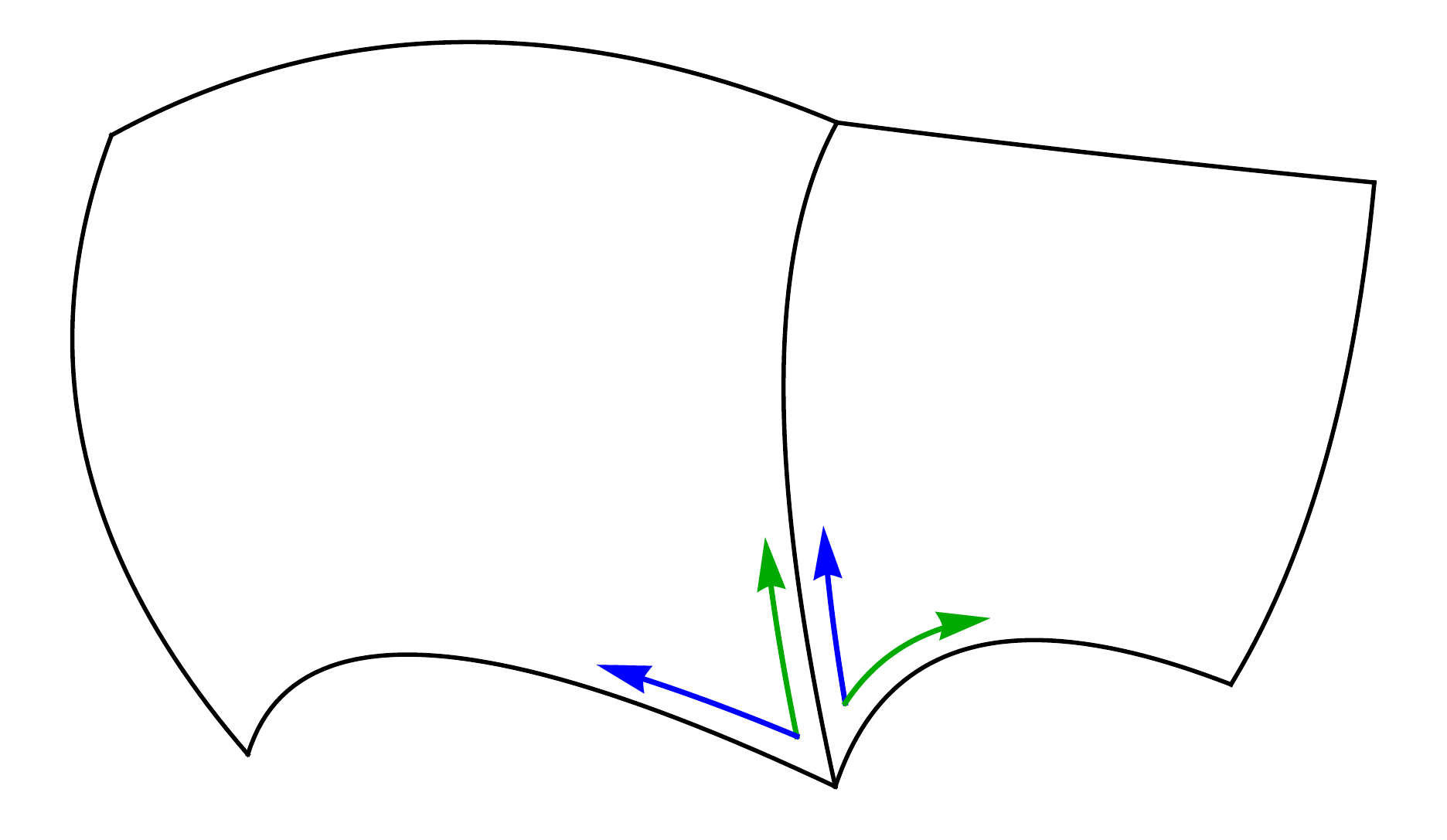}}
    \put(50,55){${\Omega}^{(i_2)}$}
    \put(140,55){${\Omega}^{(i_1)}$}
    \put(95,80){$\Sigma^{(i)}$}
  \end{picture}
  \hspace{10pt}
  \begin{picture}(160,160)
    \put(0,160){\includegraphics[width=.35\textwidth,angle=270]{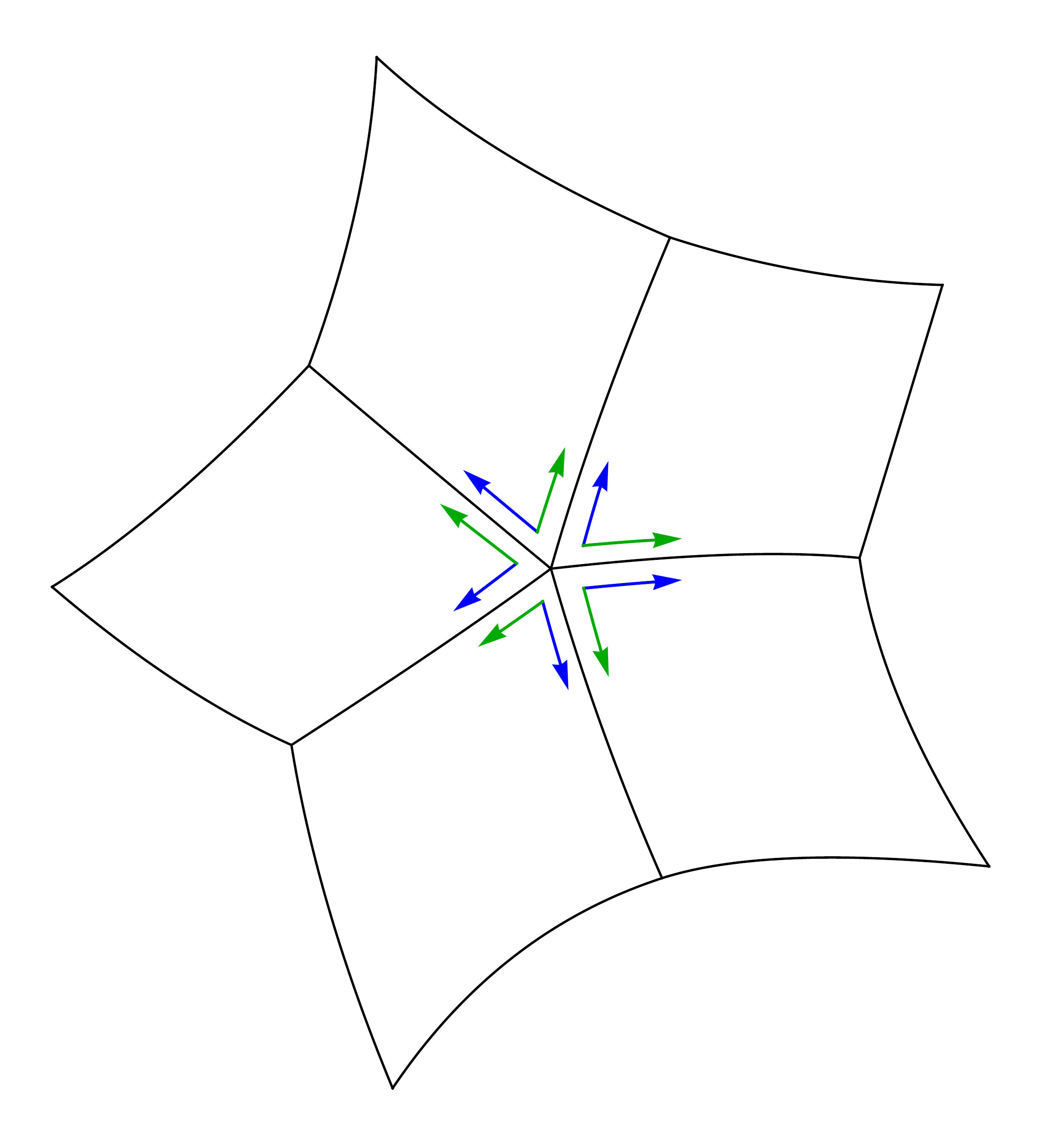}}
    \put(110,40){$\ldots$}
    \put(70,50){$\f{x}^{(i)}$}
    \put(130,75){${\Omega}^{(i_2)}$}
    \put(118,92){$\Sigma^{(i_3)}$}
    \put(85,115){${\Omega}^{(i_4)}$}
    \put(50,85){$\ldots$}
  \end{picture}
\caption{Representation in standard form~(cf.~\cite{KaSaTa17c}) of two patches~$\Omega^{(i_1)}$ and $\Omega^{(i_2)}$ with the common 
interface~$\Sigma^{(i)}$~(left) and of the patches~$\Omega^{(i_2)}$, $\Omega^{(i_4)}$, $\ldots$, $\Omega^{(i_{2\nu})}$ possessing the common vertex~$\f{x}^{(i)}$~(right).}
\label{fig:standard_representation}
\end{figure}
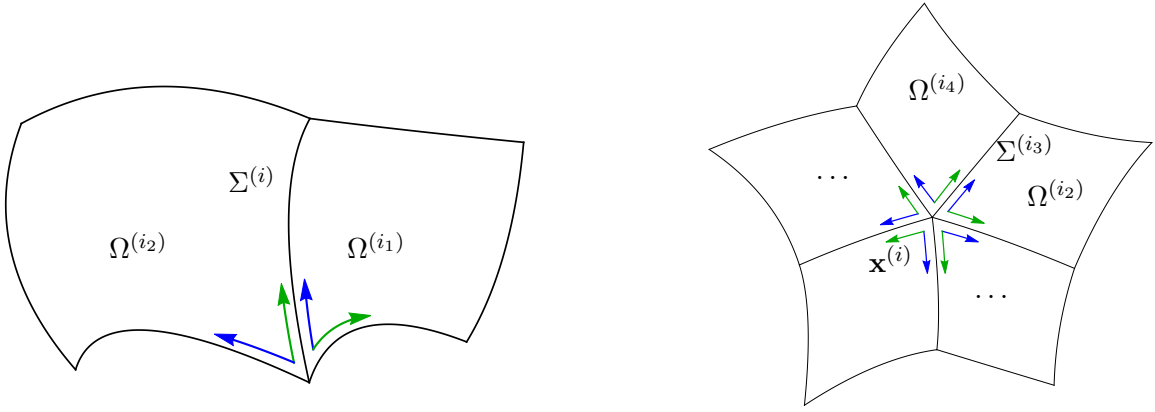

Consider an interface~$\Sigma^{(i)}$, $i \in \mathcal{I}_{\Sigma}^{\circ}$. Let $\Omega^{(i_1)}$ and $\Omega^{(i_2)}$, 
$i_1,i_2 \in \mathcal{I}_{\Omega}$, be the two neighboring patches with $\Sigma^{(i)} \subset \overline{\Omega^{(i_1)}} \cap \overline{\Omega^{(i_2)}}$. The two 
associated geometry mappings~$\f{F}^{(i_1)}$ and $\f{F}^{(i_2)}$ can be always reparameterized (if necessary) into standard form (cf.~\cite{KaSaTa17c}), which just means 
that the common interface~$\Sigma^{(i)}$ is given by
\begin{equation} \label{eq:interface_standard}
 \f{F}^{(i_1)}(0,\xi) = \f{F}^{(i_2)}(\xi,0), \mbox{ }\xi \in (0,1),
\end{equation}
see~Fig.~\ref{fig:standard_representation}~(left). 

There exist uniquely determined functions~$\alpha^{(i,i_1)} :[0,1] \rightarrow \R$, 
$\alpha^{(i,i_2)} : [0,1] \rightarrow \R$ and $\beta^{(i)} : [0,1] \rightarrow \R$ up to a common function~$\gamma^{(i)}$ 
(with $\gamma^{(i)}(\xi) \neq 0$), which are given by
\begin{equation}
   \begin{aligned} \label{eq:gluing_functions}
  	\alpha^{(i,i_1)} (\xi)  & = \gamma^{(i)} (\xi)  \det  \left [\begin{array}{ll}
          \Du \f{F}^{(i_1)}(0,\xi) & 
        \Dv \f{F}^{(i_1)}(0,\xi) 
         \end{array}\right ],\\
 	\alpha^{(i,i_2)} (\xi)  & =  \gamma^{(i)} (\xi)  \det  \left [\begin{array}{ll}
         \Du \f{F}^{(i_2)}(\xi,0) & 
        \Dv \f{F}^{(i_2)}(\xi,0) 
        \end{array}\right ],\\
 \beta^{(i)} (\xi) & =  \gamma^{(i)} (\xi) \det  \left [\begin{array}{ll}
         \Dv \f{F}^{(i_2)}(\xi,0)  & 
        \Du \f{F}^{(i_1)}(0,\xi) 
         \end{array}\right ],
   \end{aligned}
  \end{equation}
and satisfy for all $\xi \in [0,1]$
\begin{equation} \label{eq:alpha_cond}
 \alpha^{(i,i_1)}  (\xi) \alpha^{(i,i_2)}  (\xi) > 0
\end{equation}
and
\begin{equation} \label{eq:geometry_cond}
\alpha^{(i,i_1)} (\xi)  \Dv \f{F}^{(i_2)}(\xi,0)  +
        \alpha^{(i,i_2)}(\xi) \Du  \f{F}^{(i_1)}(0,\xi) + \beta^{(i)} (\xi)
        \Dv  \f{F}^{(i_1)}(0,\xi)  =\boldsymbol{0}.
\end{equation}
In addition, there exists non-unique functions~$\beta^{(i,i_1)} :[0,1] \rightarrow \R$ and $\beta^{(i,i_2)}: [0,1] \rightarrow \R$ such that
\begin{equation} \label{eq:beta_cond}
 \beta^{(i)} (\xi)= \alpha^{(i,i_1)} (\xi) \beta^{(i,i_2)}(\xi)+ \alpha^{(i,i_2)}(\xi)\beta^{(i,i_1)}(\xi),
\end{equation}
see e.g.~\cite{CoSaTa16, Pe02}. The functions~$\alpha^{(i,i_1)}$, $\alpha^{(i,i_2)}$, $\beta^{(i,i_1)}$ and $\beta^{(i,i_2)}$ are called the gluing data for the 
interface~$\Sigma^{(i)}$.

In the remainder of the paper, we will restrict ourselves to a specific class of multi-patch geometries, called analysis-suitable $G^1$ multi-patch parameterizations, 
which are needed to generate $C^1$ isogeometric spaces with optimal approximation properties, cf.~\cite{CoSaTa16}. 
\begin{definition}[Analysis-suitable $G^1$ multi-patch parameterization, cf.~\cite{CoSaTa16,KaSaTa17c}]
 A multi-patch geometry~$\f{F}$ is called analysis-suitable $G^1$ (in short, AS-$G^1$), if for every interface~$\Sigma^{(i)}$, $i \in \mathcal{I}_{\Sigma}^{\circ}$ there 
 exist linear polynomials~$\alpha^{(i,i_1)}$, $\alpha^{(i,i_2)}$, $\beta^{(i,i_1)}$ and $\beta^{(i,i_2)}$, with $\alpha^{(i,i_1)}$ and $\alpha^{(i,i_2)}$ relatively 
 prime, such that \eqref{eq:alpha_cond}, \eqref{eq:geometry_cond} and \eqref{eq:beta_cond} hold.
\end{definition}
Furthermore, for each interface~$\Sigma^{(i)}$, $i \in \mathcal{I}_{\Sigma}^{\circ}$, the linear gluing data~$\alpha^{(i,i_1)}$, $\alpha^{(i,i_2)}$, 
$\beta^{(i,i_1)}$ and $\beta^{(i,i_2)}$ is selected by minimizing the terms
\[
 ||\alpha^{(i,i_1)}-1 ||^{2}_{L_{2}([0,1])} +  ||\alpha^{(i,i_2)} -1 ||^{2}_{L_{2}([0,1])}
\]
and
\[
  ||\beta^{(i,i_1)} ||^{2}_{L_{2}([0,1])} +  ||\beta^{(i,i_2)} ||^{2}_{L_{2}([0,1])},
\]
see~\cite{KaSaTa17c}, which implies in case of parametric continuity (that is, $\beta^{(i)} \equiv 0$ and $\alpha^{(i,i_1)}=\alpha^{(i,i_2)}$) 
$\beta^{(i,i_1)} = \beta^{(i,i_2)} \equiv 0$ and $\alpha^{(i,i_1)} \equiv \alpha^{(i,i_2)} \equiv 1$.

\begin{remark}
Given a $C^1$ isogeometric function space (as defined in Section \ref{sec:C1-igfunc}) over a regular $C^0$ multi-patch parametrization. Then, for $r\leq p-2$, analysis-suitable 
$G^1$ multi-patch parameterizations are the only configurations that allow optimal approximation properties under $h$-refinement.

The reason for this is, that when the degree of the gluing data is assumed to be larger than one, there exists a configuration, such that the approximation order of function values 
and/or gradients along the interface is reduced. This is a direct consequence of Theorem 3 in \cite{CoSaTa16}. 
\end{remark}

Piecewise bilinear multi-patch parameterizations are one simple example of AS-$G^1$ multi-patch geometries~(cf. \cite{KaBuBeJu17,KaViJuBi15}), but the class of 
AS-$G^1$ multi-patch parameterizations is much wider, see e.g.~\cite{KaSaTa17b}. In Section~\ref{sec:constrASG1geometries}, we will present two possible strategies to 
construct from given non-AS-$G^1$ multi-patch geometries parameterizations which are AS-$G^1$-continuous.

\subsection{Analysis-suitable $G^1$ parameterization: construction} \label{sec:constrASG1geometries}

We describe the two approaches~\cite{KaSaTa17b,KaViJuBi15} to generate from a given initial non-AS-$G^1$ multi-patch geometry~$\widetilde{\f{F}}$ an AS-$G^1$ multi-patch 
geometry~$\f{F}$. We assume that the associated parameterizations~$\widetilde{\f{F}}^{(i)}$, $i \in \mathcal{I}_{\Omega}$, of the non-AS-$G^1$ 
geometry~$\widetilde{\f{F}}$ belong to the space $\left(\mathcal{S}_{\widetilde{h}}^{\widetilde{\f{p}},\widetilde{\f{r}}}\right)^{2}$ and are regularly parameterized. 
The goal is to construct a multi-patch geometry~$\f{F}$ consisting of 
parameterizations~$\f{F}^{(i)} \in \left( \mathcal{S}_{h}^{\f{p},\f{r}} \right)^{2}$, $i \in \mathcal{I}_{\Omega}$, with 
$\mathcal{S}_{h}^{\f{p},\f{r}} \supseteq \mathcal{S}_{\widetilde{h}}^{\widetilde{\f{p}},\widetilde{\f{r}}}$, possessing B-spline representations of the form
\begin{equation*}
 \f{F}^{(i)}(\xi_1,\xi_2) = \sum_{j_1=0}^{N-1} \sum_{j_2=0}^{N-1} \f{c}_{i_1,i_2}^{(i)} b_{(j_1,j_2)}(\xi_1,\xi_2), \mbox{ } 
 (\xi_1,\xi_2) \in [0,1]^2, 
\end{equation*}
with control points~$\f{c}_{i_1,i_2}^{(i)} \in \R^2$, such that $\f{F}$ is AS-$G^1$-continuous and that $\f{F}$ approximates $\widetilde{\f{F}}$ as good as possible.  
Below, we assume that for each edge~$\Sigma^{(i)}$, $i \in \mathcal{I}_{\Sigma}$, and for each vertex~$\f{x}^{(i)}$, $i \in \mathcal{I}_{\chi}$, the associated geometry 
mappings~$\f{F}^{(i_k)}$ are always given in standard form~\eqref{eq:interface_standard} and~\eqref{eq:vertex_standard}, respectively, compare 
also~Fig.~\ref{fig:standard_representation}. This is valid as well for the corresponding parameterizations~$\widetilde{\f{F}}^{(i_k)}$ of $\widetilde{\f{F}}$.

\subsubsection{The piecewise bilinear fitting approach~\cite{KaViJuBi15}}

Given a non-AS-$G^1$ 
multi-patch geometry~$\widetilde{\f{F}}$ with the associated parameterizations~$\widetilde{\f{F}}^{(i)}$, $i \in \mathcal{I}_{\Omega}$, we first choose a multi-patch 
geometry~$\widehat{\f{F}}$ consisting of bilinear parameterizations~$\widehat{\f{F}}^{(i)}$, $i \in \mathcal{I}_{\Omega}$, which roughly describes the initial 
geometry~$\widetilde{\f{F}}$. Then, following \cite{KaViJuBi15}, we look for a suitable approximation ~$\f{F}$ of $\widetilde{\f{F}}$, of the form  
$\f{F}= \f{u} \circ\widehat{\f{F}}$, where $\f{u}$ is a $C^1$ isogeometric vector field representing a mapping of the bilinear geometry into the final one. By the construction of Section~\ref{sec:C1ASG1space}, an 
explicit basis for $\f{u} $ is available.  Since $\widehat{\f{F}}$ is AS-$G^1$,  and  $\f{F}$ has the same gluing data by construction,  $\f{F}$  is AS-$G^1$ as required.
An example of a mapped piecewise bilinear multi-patch parameterization is given 
in~\cite[Example~6]{KaBuBeJu17} or in~\cite[Appendix~A]{CoSaTa16}, where the resulting domain is a multi-patch NURBS.

A drawback of the method is the limitation to multi-patch geometries which have to allow a rough estimation by a piecewise bilinear multi-patch parameterization. Furthermore, 
the approach cannot be used to generate AS-$G^1$ multi-patch geometries determining multi-patch domains with a smooth boundary. A more advanced technique, 
which provides amongst others the design of such multi-patch geometries, cf.~\cite[Example~3]{KaSaTa17b}, is described in the following subsection. 
 
\subsubsection{The AS-$G^1$ fitting approach~\cite{KaSaTa17b}} \label{subsec:approach_two}

This method allows the construction of an AS-$G^1$ multi-patch geometry~$\f{F}$, which interpolates the boundary, the vertices and the first derivatives at the vertices of the 
initial non-AS-$G^1$ multi-patch geometry~$\widetilde{\f{F}}$, and which is as close as possible to~$\widetilde{\f{F}}$. The construction of $\f{F}$ is divided into the following steps: 
\begin{itemize}
 \item[\emph{Step~$1$}:] For each interface~$\Sigma^{(i)}$, $i \in \mathcal{I}_{\Sigma}^{\circ}$, of the desired multi-patch geometry~$\f{F}$, 
 we precompute the gluing data of $\f{F}$ at the interface~$\Sigma^{(i)}$, that is, $\alpha^{(i,i_1)}$, $\alpha^{(i,i_2)}$, $\beta^{(i,i_1)}$ and $\beta^{(i,i_2)}$, by 
 linearizing the corresponding non-linear gluing data of $\widetilde{\f{F}}$. Let $\widetilde{\alpha}^{(i,i_1)}$, $\widetilde{\alpha}^{(i,i_1)}$ and $\widetilde{\beta}^{(i)}$ be 
 the gluing functions~\eqref{eq:gluing_functions} for the parameterizations ~$\widetilde{\f{F}}^{(i_1)}$ and $\widetilde{\f{F}}^{(i_2)}$ for $\widetilde{\gamma}^{(i)}(\xi)\equiv 1$. 
 Then, the linear functions $\alpha^{(i,i_1)}$ and $\alpha^{(i,i_2)}$ are obtained by
 \[
  \alpha^{(i,i_k)}(\xi)= \widetilde{\alpha}^{(i,i_k)}(0) \; (1-\xi) + \widetilde{\alpha}^{(i,i_k)}(1) \; \xi,  \mbox{ }k=1,2,
 \]
and the linear functions
\[
 \beta^{(i,i_k)} (\xi) = b_0^{(i,i_k)} (1-\xi) + b_{1}^{(i,i_k)} \xi , \mbox{ }b_{0}^{(i,i_k)}, b_{1}^{(i,i_k)} \in \R, \mbox{ }k=1,2,
\]
are computed by minimizing the term
\begin{equation*}
\begin{array}{l} 
 \int_{0}^{1} \|\widetilde{\beta}^{(i)} - \underbrace{(\alpha^{(i,i_1)} \beta^{(i,i_2)} +\alpha^{(i,i_2)} \beta^{(i,i_1)} )}_{\beta^{(i)}}   \|^{2} \mathrm{d}\xi \\
 + \lambda_{\beta} \left(  \int_{0}^{1} \|\beta^{((i,i_1))} \|^{2} \mathrm{d}\xi 
  + \int_{0}^{1} \|\beta^{(i,i_2)} \|^{2} \mathrm{d}\xi   \right) 
 \rightarrow \min_{(b_{0}^{(i,i_1)},b_{1}^{(i,i_1)},b_{0}^{(i,i_2)},b_{0}^{(i,i_2)})}
\end{array}
\end{equation*}
with respect to the linear constrains
\[
\beta^{(i)}(0)=\widetilde{\beta}^{(i)}(0) \mbox{ and }\beta^{(i)}(1)=\widetilde{\beta}^{(i)}(1),
\]
using a non-negative weight $\lambda_\beta$.
 \item[\emph{Step~$2$}:] We determine for the spline coefficients~$\f{c}_{j_1,j_2}^{(i)}$ of the multi-patch geometry~$\f{F}$ three different types of linear constraints, denoted by 
 $\mathcal{L}^{\circ}_{\Sigma}$, $\mathcal{L}_{\Sigma}^{\Gamma}$ and $\mathcal{L}_{\chi}$, which will be used in Step~$3$ to construct the AS-$G^1$ multi-patch geometry~$\f{F}$. 
 The constraints~$\mathcal{L}_{\Sigma}^{\circ}$ are called AS-$G^1$ constraints and will ensure that the resulting multi-patch geometry~$\f{F}$ will be AS-$G^1$-continuous. 
 For each interface~$\Sigma^{(i)}$, we require that the geometry mappings~$\f{F}^{(i_1)}$ and $\f{F}^{(i_2)}$ have to satisfy the condition~\eqref{eq:geometry_cond} for the precomputed 
 gluing data~$\alpha^{(i,i_1)}$, $\alpha^{(i,i_2)}$, $\beta^{(i,i_1)}$ and $\beta^{(i,i_2)}$ from Step~1. The so-called boundary constrains~$\mathcal{L}_{\Sigma}^{\Gamma}$ will guarantee
 that the multi-patch geometry~$\f{F}$ will coincide with   $\widetilde{\f{F}}$ on  each boundary 
 edge~$\Sigma^{(i)}$, $i \in \mathcal{I}_{\Sigma}^{\Gamma}$. 
 Finally, the so-called vertex constraints~$\mathcal{L}_{\chi}$ will ensure that $\f{F}$ will interpolate the vertices and the first derivatives of $\widetilde{\f{F}}$  at each 
 vertex~$\f{x}^{(i)}$, $i \in \mathcal{I}_{\chi}$.
All these constrains are linear, and are compatible with each other.
 \item[\emph{Step~$3$}:] Let $\f{c}$ be the vector of all control points~$\f{c}_{j_{1},j_{2}}^{(i)}$ of the multi-patch geometry~$\f{F}$. Then, the AS-$G^1$ multi-patch geometry~$\f{F}$ 
 is finally constructed by minimizing the objective function
 \begin{equation*} 
 \mathcal{F}_{2}(\f{c}) + \lambda_L \mathcal{F}_{L}(\f{c}) + 
 \lambda_U \mathcal{F}_{U}(\f{c})  \rightarrow \min_{\f{c}}
\end{equation*}
 with respect to the linear constraints~$\mathcal{L}_{\Sigma}^{\circ}$, $\mathcal{L}_{\Sigma}^{\Gamma}$ and $\mathcal{L}_{\chi}$, using non-negative weights $\lambda_{L}$ and 
 $\lambda_{U}$. While the quadratic functional $\mathcal{F}_{2}$ given by 
 \[
\mathcal{F}_{2}(\f{c}) =  \sum_{i \in \mathcal{I}_{\Omega}} \int_{[0,1]^{2}} \|\f{F}^{(i)} - \widetilde{\f{F}} \mbox{}^{(i)}\|^{2} \,\mathrm{d}\xi_1 \,\mathrm{d}\xi_2,
\]
 will ensure that the resulting geometry~$\f{F}$ approximates~$\widetilde{\f{F}}$, the so-called parametric length functional~$\mathcal{F}_{L}$ and the so-called 
 uniformity functional~$\mathcal{F}_{U}$ given by
 \[
 \mathcal{F}_{L}(\f{c}) = \sum_{i \in \mathcal{I}_{\Omega}} \int_{[0,1]^{2}} \left(\|\Du \f{F}^{(i)} \|^{2} + \|\Dv \f{F}^{(i)} \|^{2}\right)   \mathrm{d}\xi_1 \,\mathrm{d}\xi_2,
\]
and 
\[
 \mathcal{F}_{U}(\f{c}) =\sum_{i \in \mathcal{I}_{\Omega}} \int_{[0,1]^{2}} \left(\|\Du^{2} \f{F}^{(i)} \|^{2} + 2\|\Du \Dv \f{F}^{(i)} \|^{2} 
 + 
 \|\Dv^{2} \f{F}^{(i)} \|^{2}\right)  \mathrm{d}\xi_1 \,\mathrm{d}\xi_2 ,
\]
 respectively, will be needed to construct parameterizations of good quality.
 \end{itemize}
 In case that the quality of the resulting AS-$G^1$ multi-patch geometry~$\f{F}$ is not good enough, and some of the single parameterizations~$\f{F}^{(i)}$, $i \in \mathcal{I}_{\Omega}$, 
 are even singular, the use of a sufficiently refined spline space~$\mathcal{S}_{h}^{\f{p},\f{r}}$ solve, in practice,  this issue.  Two instances of constructed AS-$G^1$ multi-patch 
 geometries using this approach are given in Fig.~\ref{fig:AS_G1_geometries}.

\section{An isogeometric $C^1$ space} \label{sec:C1ASG1space}

We first introduce the concept of $C^1$ isogeometric spline spaces over (general) multi-patch geometries, and then present the construction~\cite{KaSaTa17c} to generate 
a particular $C^1$ isogeometric spline space over a given AS-$G^1$ multi-patch geometry.

\subsection{Space of $C^1$ isogeometric functions}\label{sec:C1-igfunc}

The space of $C^1$ isogeometric functions with respect to the multi-patch geometry~$\f{F}$ is defined as 
\begin{equation} \label{eq:V0}
 \mathcal{V}^{1}=\left \{ \varphi_h \in  C^1( \overline\Omega ) \, | \, \text{ for all } i \in
   \mathcal{I}_{\Omega}, \, f_h^{(i)} = \varphi_h \circ \f
   F^{(i)} \in  \mathcal{S}_{h}^{\f{p},\f{r}} \right\}.
  \end{equation}
This space can be characterized by the equivalence of the $C^1$-smoothness of an isogeometric function and the $G^1$-smoothness of its associated graph, or more 
precisely, $\varphi_h \in \mathcal{V}^{1}$ if and only if the graph of $\varphi_h$ is $G^1$-smooth at all interfaces~$\Sigma^{(i)}$, $i \in \mathcal{I}_{\Sigma}^{\circ}$, 
cf.~\cite{CoSaTa16,Pe15,KaViJuBi15}. Note that for an isogeometric function~$\varphi_h$, 
its graph~$\Phi \subset \Omega \times \R$ is the collection of the single graph surface patches
\[
 \Phi^{(i)}: [0,1]^2 \rightarrow \Omega^{(i)} \times \R, \mbox{ } \Phi^{(i)}(\xi_1,\xi_2)= \left [
    \begin{array}{c}
      \f F^{(i)}(\xi_1,\xi_2)\\ f_h^{(i)}(\xi_1,\xi_2) 
    \end{array}
\right ]  ,\mbox{ }i \in \mathcal{I}_{\Omega},
\]
with $f_h^{(i)}= \varphi_h \circ \f{F}^{(i)}$. Then, an isogeometric function $\varphi_h$ belongs to the space~$\mathcal{V}^1$, if and only if for each 
interface~$\Sigma^{(i)}$, $i \in \mathcal{I}^{\circ}_{\Sigma}$, assuming that the two associated neighboring geometry mappings~$\f{F}^{(i_1)}$ and $\f{F}^{(i_2)}$ are 
given in standard form~\eqref{eq:interface_standard}, {there exist gluing data satisfying \eqref{eq:alpha_cond} and \eqref{eq:beta_cond} such that conditions 
\eqref{eq:interface_standard} and \eqref{eq:geometry_cond} are satisfied not only for the parametrizations $\f F^{(i_1)}$, $\f F^{(i_2)}$ but also for the graph surfaces 
$\Phi^{(i_1)}$, $\Phi^{(i_2)}$. 

Since the two geometry mappings~$\f{F}^{(i_1)}$ and $\f{F}^{(i_2)}$ already uniquely determine (up to a common function~$\gamma^{(i)}$) the functions~$\alpha^{(i,i_1)}$, 
$\alpha^{(i,i_2)}$ and $\beta^{(i)}$, compare Section~\ref{subsec:analysis-suitable}, we obtain that $\varphi_h \in \mathcal{V}^{1}$ if and only if the last component of the graph 
surfaces satisfies the equations \eqref{eq:interface_standard} and \eqref{eq:geometry_cond}, that is,} for $\xi \in [0,1]$
\begin{equation*}
f_h^{(i_1)}(0,\xi_1)  = f_h^{(i_2)}(\xi,0) = (\varphi_h |_{\Sigma^{(i)}}) \circ \f F^{(i_k)}
\end{equation*}
and
\begin{equation} \label{eq:condC1_type1}
 \alpha^{(i,i_1)} (\xi)  \Dv f_h^{(i_2)}(\xi,0)  +
        \alpha^{(i,i_2)}(\xi) \Du  f_h^{(i_1)}(0,\xi) + \beta^{(i)} (\xi)
        \Dv  f_h^{(i_1)}(0,\xi)  = 0,
\end{equation}
or equivalently to~\eqref{eq:condC1_type1}
\[
 \frac{\Du  f_h^{(i_1)}(0,\xi) + \beta^{(i,i_1)}(\xi) \Dv  f_h^{(i_1)}(0,\xi)  }{\alpha^{(i,i_1)}(\xi)}  = 
 - \frac{ \Dv f_h^{(i_2)}(\xi,0) +\beta^{(i,i_2)}(\xi) \Du  f_h^{(i_2)}(\xi,0)   }{\alpha^{(i,i_2)}(\xi)} {= (\nabla \varphi_h \cdot \f d |_{\Sigma^{(i)}}) \circ \f F^{(i_k)},}
\]
where $\f d$ is a suitable vector that is not tangential to the interface, see e.g.~\cite{CoSaTa16,KaSaTa17a}. These $C^1$-conditions were used to generate $C^1$ isogeometric spline 
spaces over general analysis-suitable $G^1$ multi-patch geometries, see \cite{KaSaTa17a} for the case of two-patches and \cite{KaSaTa17c} for the multi-patch case. In the following 
subsection, we will summarize the construction~\cite{KaSaTa17c}.
  
\subsection{The Argyris isogeometric space}\label{subsec:Argyris}

We give a survey of the method~\cite{KaSaTa17c} for the design of a specific $C^1$ isogeometric spline space over a given AS-$G^1$ multi-patch geometry~$\f{F}$. 
The proposed $C^1$ space~$\mathcal{A}$ is called Argyris (quadrilateral) isogeometric space, since it possesses similar degrees-of-freedom as the classical Argyris 
triangular finite element space~\cite{argyris1968tuba}, see~\cite{KaSaTa17c} for more details. The space~$\mathcal{A}$ is a subspace of the 
entire $C^1$ isogeometric space~$\mathcal{V}^{1}$ maintaining the optimal order of approximation of the space~$\mathcal{V}^{1}$ for the traces and normal derivatives 
along the interfaces, and is much easier to investigate and to construct than the space~$\mathcal{V}^{1}$. E.g., the dimension of $\mathcal{A}$ does not depend on 
the geometry, which is in contrast to the dimension of the space~$\mathcal{V}^{1}$, cf.~\cite{KaSaTa17a} for the two-patch case. For the construction of $\mathcal{A}$, we need 
a minimal resolution within the patches given by $h \leq \frac{p-r-1}{4-r}$.

The $C^1$ isogeometric space~$\mathcal{A}$ is constructed as the direct sum of subspaces referring to the single patch-interior, edge and vertex components, that is,
\[
 \mathcal{A} = \left( \bigoplus_{i \in \mathcal{I}_{\Omega}} \mathcal{A}^{\circ}_{\Omega^{(i)}} \right) \oplus 
 \left( \bigoplus_{i \in \mathcal{I}_{\Sigma}} \mathcal{A}^{\circ}_{\Sigma^{(i)}} \right) \oplus
 \left( \bigoplus_{i \in \mathcal{I}_{\chi}} \mathcal{A}_{\f{x}^{(i)}}  \right) .
\]
The spaces~$\mathcal{A}^{\circ}_{\Omega^{(i)}}$, $\mathcal{A}^{\circ}_{\Sigma^{(i)}}$ and $\mathcal{A}_{\f{x}^{(i)}}$ are called patch-interior, edge and 
vertex function space, respectively. {The patch-interior functions are completely supported within one patch, the edge functions have support along the edge and are restricted to 
two patches, whereas the vertex functions have support in a neighborhood of the vertex. The different types of functions} are defined as follows: 

\paragraph{Patch-interior function space~$\mathcal{A}^{\circ}_{\Omega^{(i)}}$} Let $i \in \mathcal{I}_{\Omega}$. The space~$\mathcal{A}^{\circ}_{\Omega^{(i)}}$ is 
given as
\[
 \mathcal{A}^{\circ}_{\Omega^{(i)}} = \Span \{\mathrm{B}_{\f{j}}^{(i)}: \; \f{j} \in \{2,\ldots,N-3 \}^2  \}
\]
with
\begin{equation*}
    \mathrm{B}_{\f{j}}^{(i)} (\f{x}) =  
    \begin{cases}
 \left( b_{\f{j}} \circ \right(\f{F}^{(i)}\left)^{-1}\right)(\f{x}) & \mbox{if }\f \, \f{x} \in \overline{\Omega^{(i)}} ,\\
0 & \mbox{otherwise}.
\end{cases}
\end{equation*}
The patch-interior functions~$\mathrm{B}_{\f{j}}^{(i)}$, $\f{j} \in \{2,\ldots,N-3 \}^2$, are the ``standard'' isogeometric function with a support entirely contained 
in~$\Omega^{(i)}$ and have vanishing function values and vanishing gradients at the patch boundary~$\partial \Omega^{(i)}$. This directly implies that 
$\mathrm{B}_{\f{j}}^{(i)} \in C^1(\Omega)$ for $\f{j} \in \{2,\ldots,N-3 \}^2$. The dimension of the space~$\mathcal{A}^{\circ}_{\Omega^{(i)}}$, 
$i \in \mathcal{I}_{\Omega}$, is given by
\[
 \dim (\mathcal{A}^{\circ}_{\Omega^{(i)}}) = ((p-r)(n-1)+p-3)^2.
\]

\paragraph{Edge function space~$\mathcal{A}^{\circ}_{\Sigma^{(i)}}$}

We consider first the case of an interface~$\Sigma^{(i)}$, which means that $i \in \mathcal{I}^{\circ}_{\Sigma}$, and assume without loss of generality that the two associated 
neighboring geometry mappings~$\f{F}^{(i_1)}$ and $\f{F}^{(i_2)}$ are given in standard form~\eqref{eq:interface_standard}. Let $b_{j}^{+}$, $j=0,\ldots, N_0 -1$, with 
$N_0=p+(n-1)(p-r-1)+1$, be the B-splines of the univariate spline space~$\mathcal{S}^{p,r+1}_{h}$, and let $b_{j}^{-}$, $j=0,\ldots,N_1-1$, with 
$N_1=p+(n-1)(p-r-1)$, be the B-splines of the univariate spline space~$\mathcal{S}^{p-1,r}_h$. The space~$\mathcal{A}^{\circ}_{\Sigma^{(i)}}$ is defined as
\[
 \mathcal{A}^{\circ}_{\Sigma^{(i)}} = \Span \{ \overline{\mathrm{B}}_{(j_1,j_2)}^{(i)}: \; j_1=0,\ldots, N_{j_2}-1,\, j_2=0,1\}
\]
with
\begin{equation*}
    \overline{\mathrm{B}}_{(j_1,j_2)}^{(i)} (\f{x}) =  
    \begin{cases}
 \left( \overline{f}_{(j_1,j_2)}^{(i,k)} \circ \right(\f{F}^{(k)}\left)^{-1}\right)(\f{x}) & \mbox{if }\f \, \f{x} \in \overline{\Omega^{(k)}} , \mbox{ }k=i_1,i_2,\\
0 & \mbox{otherwise},
\end{cases}
\end{equation*}
where
\begin{equation}
  \begin{aligned}
 \overline{f}_{(j_1,0)}^{(i,i_1)}(\xi_1,\xi_2)  & = b_{j_1}^{+}(\xi_2)(b_0(\xi_1)+b_{1}(\xi_1)) - \beta^{(i,i_1)}(\xi_2)(b_{j_1}^{+})'(\xi_2)\frac{h}{p}b_1(\xi_1), \\
 \overline{f}_{(j_1,0)}^{(i,i_2)}(\xi_1,\xi_2)  & = b_{j_1}^{+}(\xi_1)(b_0(\xi_2)+b_{1}(\xi_2)) - \beta^{(i,i_2)}(\xi_1)(b_{j_1}^{+})'(\xi_1)\frac{h}{p}b_1(\xi_2),
  \end{aligned}\label{eq:edge-function-0}
\end{equation}
and
\begin{equation}
 \begin{aligned}
  \overline{f}_{(j_1,1)}^{(i,i_1)}(\xi_1,\xi_2) & =  \alpha^{(i,i_1)}(\xi_2)b_{j_1}^{-}(\xi_2)b_1(\xi_1), \\
  \overline{f}_{(j_1,1)}^{(i,i_2)}(\xi_1,\xi_2) & =  -\alpha^{(i,i_2)}(\xi_1)b_{j_1}^{-}(\xi_1)b_1(\xi_2).
 \end{aligned}\label{eq:edge-function-1}
\end{equation}
In case of a boundary edge~$\Sigma^{(i)}$, $i \in \mathcal{I}^{\Gamma}_{\Sigma}$, the space~$\mathcal{A}^{\circ}_{\Sigma^{(i)}}$ can be defined in a similar way. 
Assume that the associated geometry mapping is $\f{F}^{(i_1)}$, and is parameterized as in~\eqref{eq:interface_standard} for the case of two patches. The gluing 
data~$\alpha^{(i,i_1)}$ and $\beta^{(i,i_1)}$ can be simplified to $\alpha^{(i,i_1)}(\xi)=1$ and $\beta^{(i,i_1)}(\xi)=0$, which leads to functions
\begin{equation*}
    \overline{\mathrm{B}}_{(j_1,j_2)}^{(i)} (\f{x}) =  
    \begin{cases}
 \left( \overline{f}_{(j_1,j_2)}^{(i,i_1)} \circ \right(\f{F}^{(i_1)}\left)^{-1}\right)(\f{x}) & \mbox{if }\f \, \f{x} \in \overline{\Omega^{(i_1)}} ,\\
0 & \mbox{otherwise},
\end{cases}
\end{equation*}
with 
\begin{equation*}
  \begin{aligned}
 \overline{f}_{(j_1,0)}^{(i,i_1)}(\xi_1,\xi_2)  & = b_{j_1}^{+}(\xi_2)(b_0(\xi_1)+b_{1}(\xi_1)), 
  \end{aligned}
\end{equation*}
and
\begin{equation*}
 \begin{aligned}
  \overline{f}_{(j_1,1)}^{(i,i_1)}(\xi_1,\xi_2) & =  b_{j_1}^{-}(\xi_2)b_1(\xi_1). \\
 \end{aligned}
\end{equation*}
The edge functions~$\overline{\mathrm{B}}_{(j_1,j_2)}^{(i)}$, $j_1=0,\ldots, N_{j_2}-1, j_2=0,1$, are constructed in such a way that they are $C^1$-smooth across the 
interface~$\Sigma^{(i)}$, and that they span function values and cross derivative values along the edge, see \cite{KaSaTa17c} for details. In addition, the edge 
functions~$\overline{\mathrm{B}}_{(j_1,j_2)}^{(i)}$ possess a support, which is entirely contained in $\overline{\Omega^{(i_1)}} \cup \overline{\Omega^{(i_2)}}$ 
($\overline{\Omega^{(i_1)}}$ if $\Sigma^{(i)}$ is a boundary edge) in an h-dependent neighborhood of~$\Sigma^{(i)}$, and have vanishing derivatives up to second order 
at the endpoints (vertices) of the edge. This implies that $\overline{\mathrm{B}}_{(j_1,j_2)}^{(i)} \in C^{1}(\Omega)$, $j_1=0,\ldots, N_{j_2}-1, j_2=0,1$.  
The dimension of the space~$\mathcal{A}_{\Sigma^{(i)}}^{\circ}$, $i \in \mathcal{I}_{\Sigma}$, is given by
\[
 \dim (\mathcal{A}^{\circ}_{\Sigma^{(i)}}) = 2(p-r-1)(n-1)+p-9.
\]

\paragraph{Vertex function space~$\mathcal{A}_{\f{x}^{(i)}}$}

We consider a vertex~$\f{x}^{(i)}$, $i \in \mathcal{I}_{\chi}$, and denote by $\nu$ the patch valence of the vertex~$\f{x}^{(i)}$. Let $\Sigma^{(i_1)}$, $\Omega^{(i_2)}$, 
$\Sigma^{(i_3)}$, $\ldots$, $\Omega^{(i_{2\nu })}$, $\Sigma^{(i_{2\nu +1})}$, be the sequence of interfaces and patches around the vertex~$\f{x}^{(i)}$ 
in counterclockwise order, assuming in case of an inner vertex~$\f{x}^{(i)}$, $i \in \mathcal{I}_{\chi}^{\circ}$, that $\Sigma^{(i_1)}=\Sigma^{(i_{2\nu +1})}$. 
The associated geometry mappings~$\f{F}^{(i_2)}$, $\f{F}^{(i_4)}$, $\ldots$, $\f{F}^{(i_{2\nu})}$ containing the vertex~$\f{x}^{(i)}$ can be always reparameterized 
(if necessary) into standard form (cf.~\cite{KaSaTa17c}), just meaning that we have
\begin{equation} \label{eq:vertex_standard}
 \f{F}^{(i_{2k})}(0,\xi) = \f{F}^{(i_{2k+2})}(\xi,0) , \mbox{ }\xi \in [0,1 ],
\end{equation}
for $k \in \{1, \ldots, \nu -1 \}$, and additionally
\begin{equation*}
  \f{F}^{(i_{2\nu})}(0,\xi) = \f{F}^{(i_{2})}(\xi,0) , \mbox{ }\xi \in [0,1 ],
\end{equation*}
in case of an inner vertex~$\f{x}^{(i)}$, see Fig.~\ref{fig:standard_representation}~(right). This implies that
\[
 \f{x}^{(i)} = \f{F}^{(i_2)}(0,0) = \f{F}^{(i_4)}(0,0) = \ldots = \f{F}^{(i_{2\nu})}(0,0).
\]
Considering a boundary vertex~$\f{x}^{(i)}$, $i  \in \mathcal{I}_{\chi}^{\Gamma}$, we assume that the edges~$\Sigma^{(i_1)}$ and $\Sigma^{(2i_{\nu + 1})}$ are the two boundary 
edges, for which the gluing data~$\alpha^{(i_1,i_2)}$, $\alpha^{(i_{2 \nu +1},i_{2\nu})}$ and $\beta^{(i_1,i_2)}$, $\beta^{(i_{2 \nu +1},i_{2\nu})}$ 
can be simplified to $\alpha^{(i_1,i_2)}(\xi)= \alpha^{(i_{2 \nu +1},i_{2\nu})}(\xi)=1$ and $\beta^{(i_1,i_2)}(\xi) =\beta^{(i_{2 \nu +1},i_{2\nu})}(\xi)=0$. 

Before defining the space~$\mathcal{A}_{\f{x}^{(i)}}$, $i \in \mathcal{I}_{\chi}$, we need further tools  and definitions. We consider 
the basis transformations~$\{b_0, b_1 \}$ to $\{c_0,c_1 \}$, $\{b_{0}^{+}, b_{1}^{+},b_{2}^{+} \}$ to $\{c_{0}^{+}, c_{1}^{+},c_{2}^{+} \}$ and from 
$\{b_{0}^{-},b_{1}^{-}\}$ to $\{c_{0}^{-},c_{1}^{-}\}$, with
\begin{equation*}
 \partial^j_{\xi} c_i(0)=\delta_i^j \quad\mbox{for }j=0,1, \quad \partial^j_{\xi}c_{i}^{+}(0)=\delta_i^j \quad\mbox{for }j=0,\ldots,2,\quad\mbox{and}\quad   
 \partial^j_{\xi} c_i^-(0)=\delta_i^j \quad\mbox{for }j=0,1,
\end{equation*}
where $\delta_i^j$ is the Kronecker delta. For each edge~$\Sigma^{(i_k)}$, $k \in \{1,3,\ldots, 2 \nu +1\}$, we use the abbreviated notations 
\[
 \f{t}^{(i_k)} = \Dv \f{F}^{i_{k-1}}(0,\xi) = \Du \f{F}^{i_{k+1}}(\xi,0)
\]
and
\begin{equation*}
\begin{array}{lll}
  \f{d}^{(i_k)}(\xi) &=& \frac{1}{\alpha^{(i_k,i_{k-1})}(\xi)}\left( \Du \f F^{(i_{k-1})}(0,\xi) + \beta^{(i_k,i_{k-1})}(\xi) \, \Dv \f F^{(i_{k-1})}(0,\xi)\right) \\
  &=& -\frac{1}{\alpha^{(i_k,i_{k+1})}(\xi)}\left( \Dv \f F^{(i_{k+1})}(\xi,0) + \beta^{(i_k,i_{k+1})}(\xi) \, \Du \f F^{(i_{k+1})}(\xi,0)\right).
\end{array}
\end{equation*}
Note that in case of a boundary edge in each case one term does not exist. Given the vector~$\Phi = (\phi_{0,0},\phi_{1,0},\phi_{0,1},\phi_{2,0},\phi_{1,1},\phi_{0,2})$, which 
describes the $C^2$ interpolation data of an isogeometric function~$\phi$ at the vertex~$\f{x}^{(i)}$ determined by the function value~$\phi(\f{x}^{(i)}) = \phi_{0,0}$, 
the gradient $\nabla \phi = (\phi_{1,0},\phi_{0,1})$ and the Hessian
\begin{equation*}
H\phi = \left(
\begin{array}{ll}
  \phi_{2,0} & \phi_{1,1} \\
  \phi_{1,1} & \phi_{0,2} 
\end{array}\right),
\end{equation*} 
we define for each patch~$\Omega^{(i_k)}$, $k \in \{2,4 \ldots, 2\nu\}$, the functions
\begin{equation*}
 \begin{aligned}
 \ev{f}_{\Phi}^{(i_{k+1},i_{k})}(\xi_1,\xi_2) & = \sum_{j=0}^{2} d_{0,j}^{(i_{k+1},i_k)} \left(c_{j}^{+}(\xi_2)c_0(\xi_1) - 
 \beta^{(i_{k+1},i_k)}(\xi_2)(c_{j}^{+})'(\xi_2)c_1(\xi_1) \right) \\
 & + \sum_{j=0}^{1} d_{1,j}^{(i_{k+1},i_k)} \alpha^{(i_{k+1},i_{k})}(\xi_2) c_{j}^{-}(\xi_2)c_{1}(\xi_1),
 \end{aligned}
 \end{equation*}
 \begin{equation*}
 \begin{aligned}
  \ev{f}_{\Phi}^{(i_{k-1},i_{k})}(\xi_1,\xi_2) & = \sum_{j=0}^{2} d_{0,j}^{(i_{k-1},i_k)} \left(c_{j}^{+}(\xi_1)c_0(\xi_2) - 
  \beta^{(i_{k-1},i_k)}(\xi_1)(c_{j}^{+})'(\xi_1)c_1(\xi_2) \right) \\
 & - \sum_{j=0}^{1} d^{(i_{k-1},i_k)}_{1,j} \alpha^{(i_{k-1},i_{k})}(\xi_1) c_{j}^{-}(\xi_1)c_{1}(\xi_2),
 \end{aligned}
\end{equation*}
and
\begin{equation*}
 \ev{f}_{\Phi}^{(i_k)}(\xi_1,\xi_2) = \sum_{j_1=0}^{1} \sum_{j_2=0}^{1} {d}^{(i_k)}_{j_1,j_2} c_{1}(\xi_1)c_{2}(\xi_2),
\end{equation*}
with
\begin{equation*}
 d_{0,0}^{(i_{\ell},i_k)} = \phi_{0,0}, \mbox{ } d_{0,1}^{(i_{\ell},i_k)} = \nabla \phi \; \f t^{(i_{\ell})}(0), 
 \mbox{ }d_{0,2}^{(i_{\ell},i_k)} = (\f{t}^{(i_{\ell})}(0))^T  \; H\phi \; \f{t}^{(i_{\ell})}(0) + \nabla\phi \; (\f t^{(i_{\ell})})'(0), 
\end{equation*}
\begin{equation*}
  d_{1,0}^{(i_{\ell},i_k)} = \nabla \phi \; \f{d}^{(i_{\ell})}(0) , \mbox{ }
   d_{1,1}^{(i_{\ell},i_{k})} = (\f{t}^{(i_{\ell})}(0))^T \; H\phi \; \f{d}^{(i_{\ell})}(0) + \nabla\phi \;(\f{d}^{(i_{\ell})})'(0), 
\end{equation*}
for $\ell=k-1,k+1$, and
\begin{equation*}
 d_{0,0}^{(i_k)} = \phi_{0,0},\mbox{ } d_{1,0}^{(i_k)}= \nabla \phi \; \f t^{(i_{k-1})}(0),\mbox{ } d_{0,1}^{(i_k)} =  \nabla \phi \; \f t^{(i_{k+1})}(0) 
\end{equation*}
\begin{equation*}
d_{1,1}^{(i_k)} = (\f{t}^{(i_{k-1})}(0))^T  \; H\phi \; \f{t}^{(i_{k+1})}(0) + \nabla \phi \; \Du \Dv \f{F}^{(i_k)}(0,0).
\end{equation*}
Let
\begin{equation*}
 \Phi_{(0,0)} = (1,0,0,0,0,0), \mbox{ }\Phi_{(1,0)} = (0,1,0,0,0,0), \mbox{ } \Phi_{(0,1)} = (0,0,1,0,0,0),
\end{equation*}
and
\begin{equation*}
 \Phi_{(2,0)} = (0,0,0,1,0,0), \mbox{ }\Phi_{(1,1)} = (0,0,0,0,1,0), \mbox{ } \Phi_{(0,2)} = (0,0,0,0,0,1),
\end{equation*}
then the space~$\mathcal{A}_{\f{x}^{(i)}}$ is defined as
\[
 \mathcal{A}_{\f{x}^{(i)}} = \Span \{ \ev{\mathrm{B}}_{(j_1,j_2)}^{(i)}: \; 0 \leq j_1,j_2 \leq 2, \, j_1+j_2 \leq2\}
\]
with
\begin{equation*}
    \ev{\mathrm{B}}_{(j_1,j_2)}^{(i)} (\f{x}) =  
    \begin{cases}
 \left( \left( \ev{f}_{\Phi_{(j_1,j_2)}}^{(i_{k-1},i_k)} + \ev{f}_{\Phi_{(j_1,j_2)}}^{(i_{k+1},i_k)} -\ev{f}_{\Phi_{(j_1,j_2)}}^{(i_k)}   \right)
 \circ \right(\f{F}^{(i_k)}\left)^{-1}\right)(\f{x}) & \mbox{if }\f \, \f{x} \in \overline{\Omega^{(i_k)}} , \mbox{ }k=2,4,\ldots,2\nu,\\
0 & \mbox{otherwise},
\end{cases}
\end{equation*}
where the factor
\[
  \sigma = \left (  \frac{h}{ p \, \nu} \sum_{\ell = 1}^\nu \|\nabla \f{F}^{(i_{2\ell})} (0,0) \| \right ) ^{-1}
\]
is used to uniformly scale the functions with respect to the $L^{\infty}$-norm. The vertex functions~$\ev{\mathrm{B}}_{(j_1,j_2)}^{(i)}$, $0 \leq j_1,j_2 \leq 2$, $j_1+j_2$, 
are constructed in such a way that they are $C^1$-smooth across all interfaces~$\Sigma^{(i_{2k+1)}}$, $k=0,\ldots,\nu$, and that they span the function value and all derivatives 
up to second order at the vertex~$\f{x}^{(i)}$, see~\cite{KaSaTa17c} for details. Furthermore, the support of an vertex function~$\ev{\mathrm{B}}_{(j_1,j_2)}^{(i)}$ is 
entirely contained in $\cup_{k=1}^{\nu} \overline{\Omega^{i_{2k}}}$ in an $h$-dependent neighborhood of the vertex~$\f{x}^{(i)}$. Therefore, we obtain that 
$\ev{\mathrm{B}}_{(j_1,j_2)}^{(i)} \in C^{1}(\Omega)$ and $\ev{\mathrm{B}}_{(j_1,j_2)}^{(i)} \in C^{2}(\f{x}^{(i)})$ for $0 \leq j_1,j_2 \leq 2$, $j_1+j_2$. 
The dimension of the space~$\mathcal{A}_{\f{x}^{(i)}}$, $i \in \mathcal{I}_{\chi}$, is equal to
\[
\dim (\mathcal{A}_{\f{x}^{(i)}}) = 6.
\]
As already mentioned above, the dimension of the entire space~$\mathcal{A}$ does not depend on the geometry and is finally given by the sum of dimensions of all 
subspaces~$\mathcal{A}_{\Omega^{(i)}}$, $\mathcal{A}_{\Sigma^{(i)}}$ and $\mathcal{A}_{\f{x}^{(i)}}$, that is,
\[
 \dim(\mathcal{A}) = |\mathcal{I}_{\Omega}| \cdot \left((p-r)(n-1)+p-3\right)^2 + |\mathcal{I}_{\Sigma}| \cdot \left(2(p-r-1)(n-1)+p-9\right) + 
 |\mathcal{I}_{\chi}| \cdot 6.
\]
Moreover, all constructed patch-interior functions $\mathrm{B}_{(j_1,j_2)}^{(i)}$, edge functions $\overline{\mathrm{B}}_{(j_1,j_2)}^{(i)}$ and 
vertex functions $\ev{\mathrm{B}}_{(j_1,j_2)}^{(i)}$ from above form a basis of the space~$\mathcal{A}$. This is a direct consequence of the definition of the 
single functions and their supports.

\section{Beyond analysis-suitable $G^1$ parameterizations}
\label{sec:generalization}
It is possible to extend the $C^1$ basis construction to domains that are not analysis-suitable $G^1$. This is done by locally increasing the degree. This approach was employed for 
interfaces in the multi-patch framework in \cite{ChAnRa18}, based on the findings in \cite{CoSaTa16}. More extensive research was done in \cite{Pe15-2,KaPe15A,KaPe18,NgKaPe15} for 
unstructured quadrilateral meshes, where higher degree elements where used in local regions around extraordinary vertices.

\begin{figure}[ht!]
\centering
\begin{subfigure}[t]{0.25\textwidth}
\includegraphics[width=.8\textwidth]{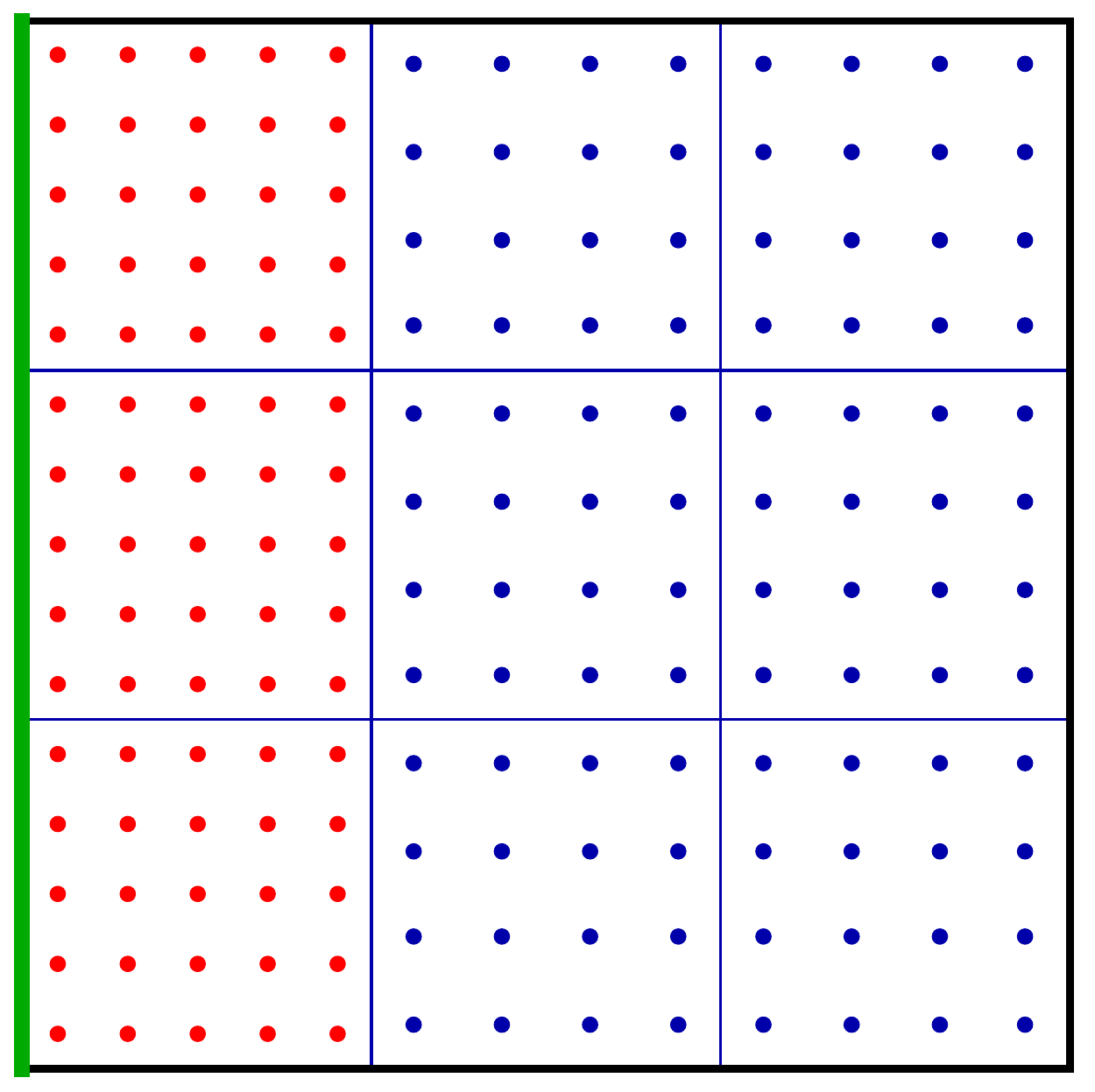}
\caption{One interface without special vertex construction.}\label{subfig:1if}
\end{subfigure}\qquad
\begin{subfigure}[t]{0.25\textwidth}
\includegraphics[width=.8\textwidth]{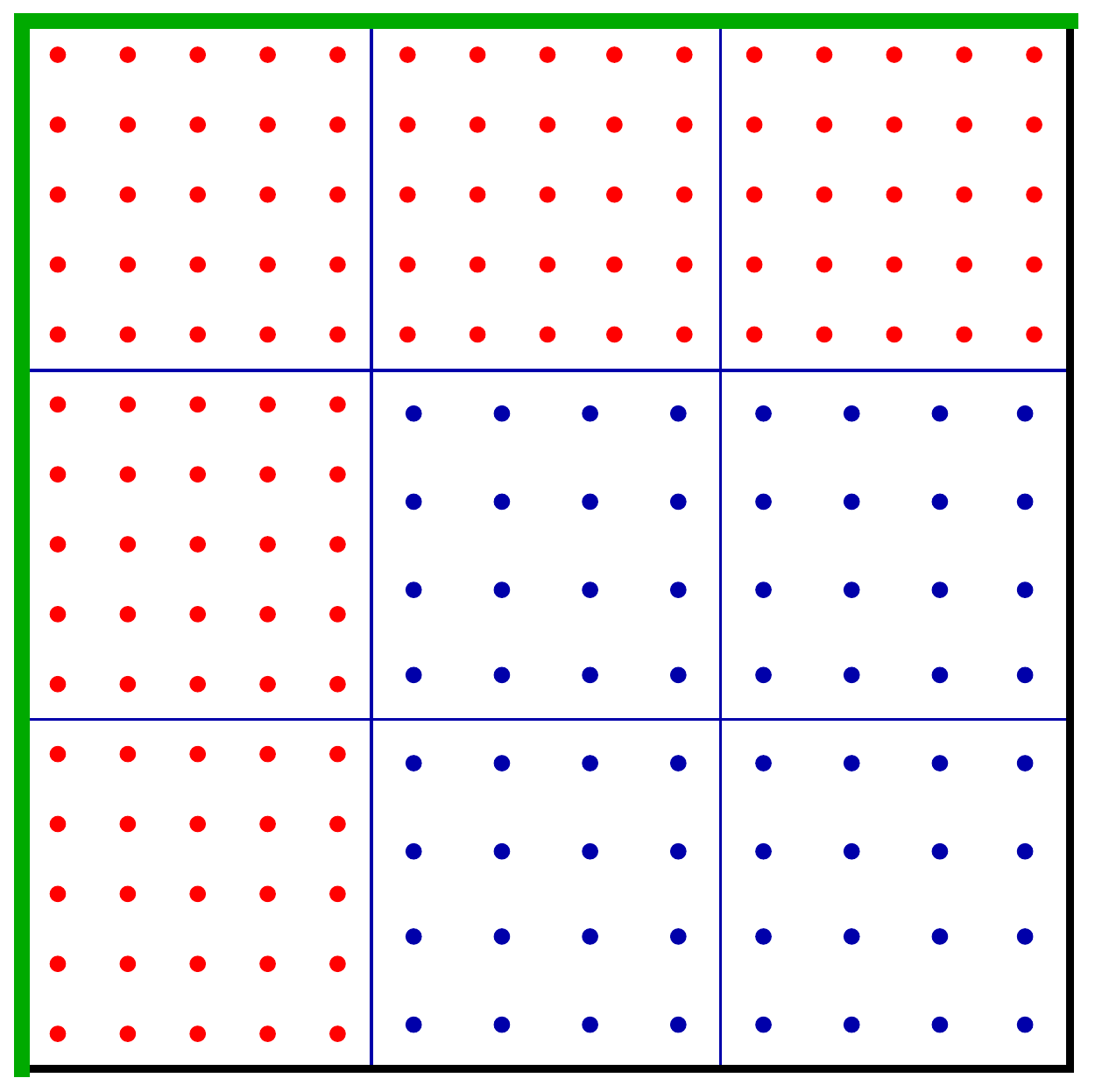}
\caption{Two interfaces; coupled at the top-left vertex.}\label{subfig:2if}
\end{subfigure}\qquad
\begin{subfigure}[t]{0.25\textwidth}
\includegraphics[width=.8\textwidth]{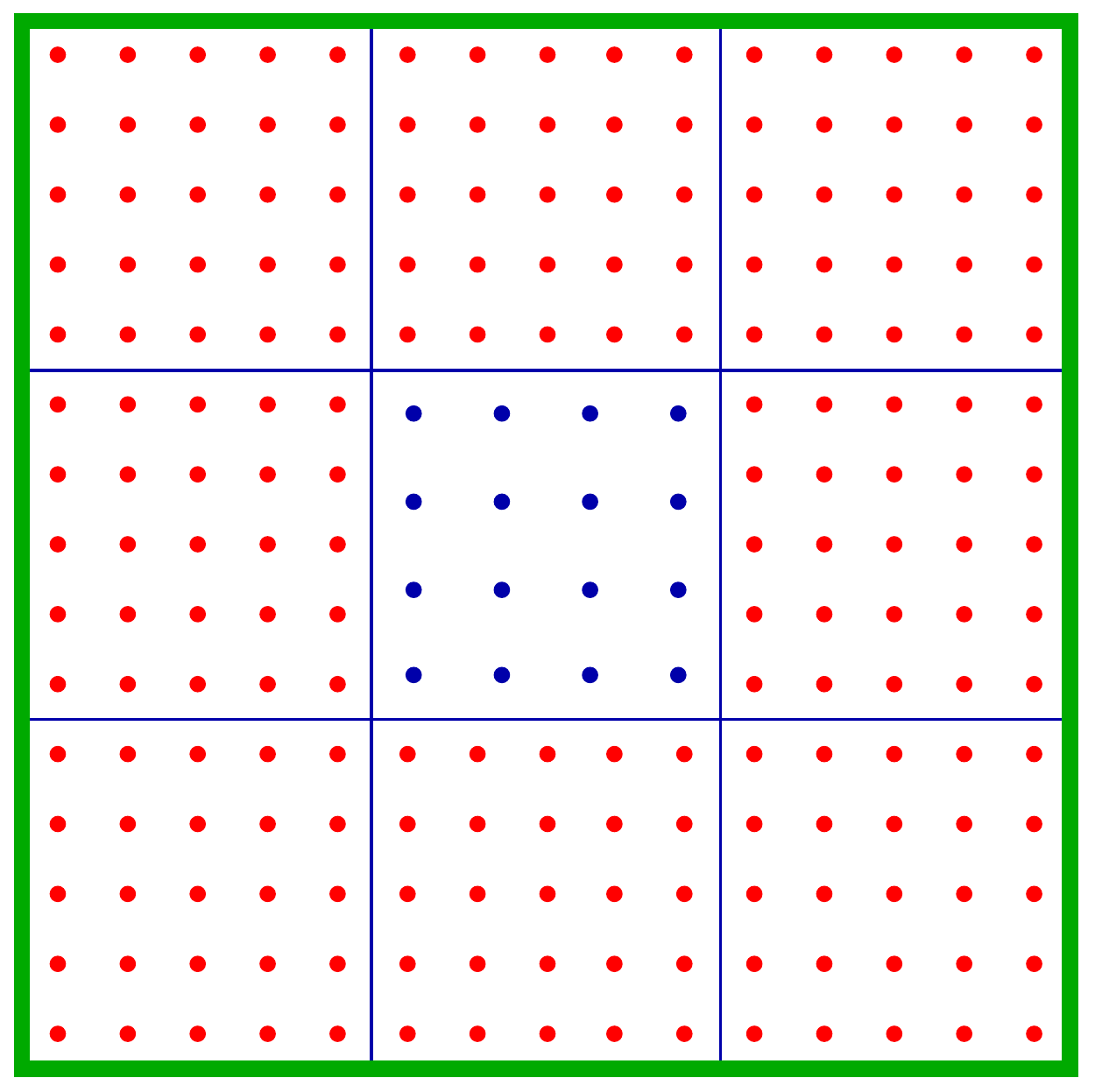}
\caption{Four interfaces; coupled at all vertices.}\label{subfig:4if}
\end{subfigure}
\caption{Spline-like patches of mixed degree (coefficients for $p=4$ in red, for $p=3$ in blue). The $G^1$ interfaces are depicted in green.}
\label{fig:multi-degree-patch}
\end{figure}

The mixed degree construction is based on the following observation. The definitions of the edge basis functions in \eqref{eq:edge-function-0} and \eqref{eq:edge-function-1} are not 
confined to the gluing data being linear functions. However, we have the following lemma.
\begin{lemma}\label{lem:edge-basis-mixed-degree}
 Given an interface $\Sigma^{(i)}$ between patches $\Omega^{(i_1)}$ and $\Omega^{(i_2)}$. Let $\alpha^{(i,i_1)},\beta^{(i,i_1)} \in \mathcal{S}^{p^*,r^*}_h$ and $b^*_0$, $b^*_{1}$ be 
 the first two basis functions in $\mathcal{S}^{p^+,r^-}_h$, where $p^+ = \max(p,p+p^*-1)$ and $ r^- = \min(r,r^*)$. Then the edge functions 
 \begin{equation*}
 \overline{f}_{(j_1,0)}^{(i,i_1)}(\xi_1,\xi_2)  = b_{j_1}^{+}(\xi_2)(b^*_0(\xi_1)+b^*_{1}(\xi_1)) - \beta^{(i,i_1)}(\xi_2)(b_{j_1}^{+})'(\xi_2)\frac{h}{p}b^*_1(\xi_1)
\end{equation*}
and
\begin{equation*}
  \overline{f}_{(j_1,1)}^{(i,i_1)}(\xi_1,\xi_2) =  \alpha^{(i,i_1)}(\xi_2)b_{j_1}^{-}(\xi_2)b^*_1(\xi_1)
\end{equation*}
 satisfy
 \[
  \overline{f}_{(j_1,0)}^{(i,i_1)},\overline{f}_{(j_1,1)}^{(i,i_1)} \in \mathcal{S}^{\f p^+,\f r^-}_h,
 \]
with $\f p^+ = (p^+,p^+)$ and $\f r^- = (r^-,r^-)$. Analogously, we have the same for the patch $\Omega^{(i_2)}$.
\end{lemma}
A simple consequence is, that if there exists quadratic gluing data, then all edge functions are in the space $\mathcal{S}^{\f p + \f 1,\f r}_h$. This significantly increases 
the flexibility of the multi-patch geometry. Note that the biquadratic G-splines defined in \cite{Re95} possess quadratic gluing data. Configurations of suitable spline-like patches 
of mixed degree are given in Figure \ref{fig:multi-degree-patch}. The blue dots signify B\'ezier coefficients of degree $3$, whereas the red dots correspond to B\'ezier coefficients 
of degree $4$. The blue line are mesh lines where parametric continuity of order $C^1$ or $C^2$ is described whereas the green edges at the boundary are of $G^1$ smoothness across 
patches and the black edges are either boundary edges, or edges where the patch can be extended with parametric continuity (at least $C^1$). Note that one may prescribe different 
continuity for different regions of the patch, e.g., only $C^1$ close to the interfaces and $C^2$ in the interior.

In this configuration of mixed degree $3$ and $4$, the coefficients corresponding to the inner elements (not neighboring the interfaces) do not influence the function value or value 
of first derivatives at the interfaces. Hence, the edge and vertex functions are completely determined by the B\'ezier elements of degree $4$.

Considering a configuration as in Figure \ref{subfig:1if} containing only one interface, the edge basis as presented in Lemma \ref{lem:edge-basis-mixed-degree} together with the 
patch interior basis obtained by resolving the $C^k$ conditions in the interior give a complete basis of the $C^1$ smooth isogeometric function space. 

When given configurations as in Figures \ref{subfig:2if} or \ref{subfig:4if} additional vertex functions can be defined by interpolation of $C^2$ data. The procedure is similar to 
the construction presented in Subsection \ref{subsec:Argyris}. In all configurations, the patch interior basis is a tensor-product of suitable univariate basis functions.

\begin{figure}[ht!]
\centering
\includegraphics[width=.3\textwidth]{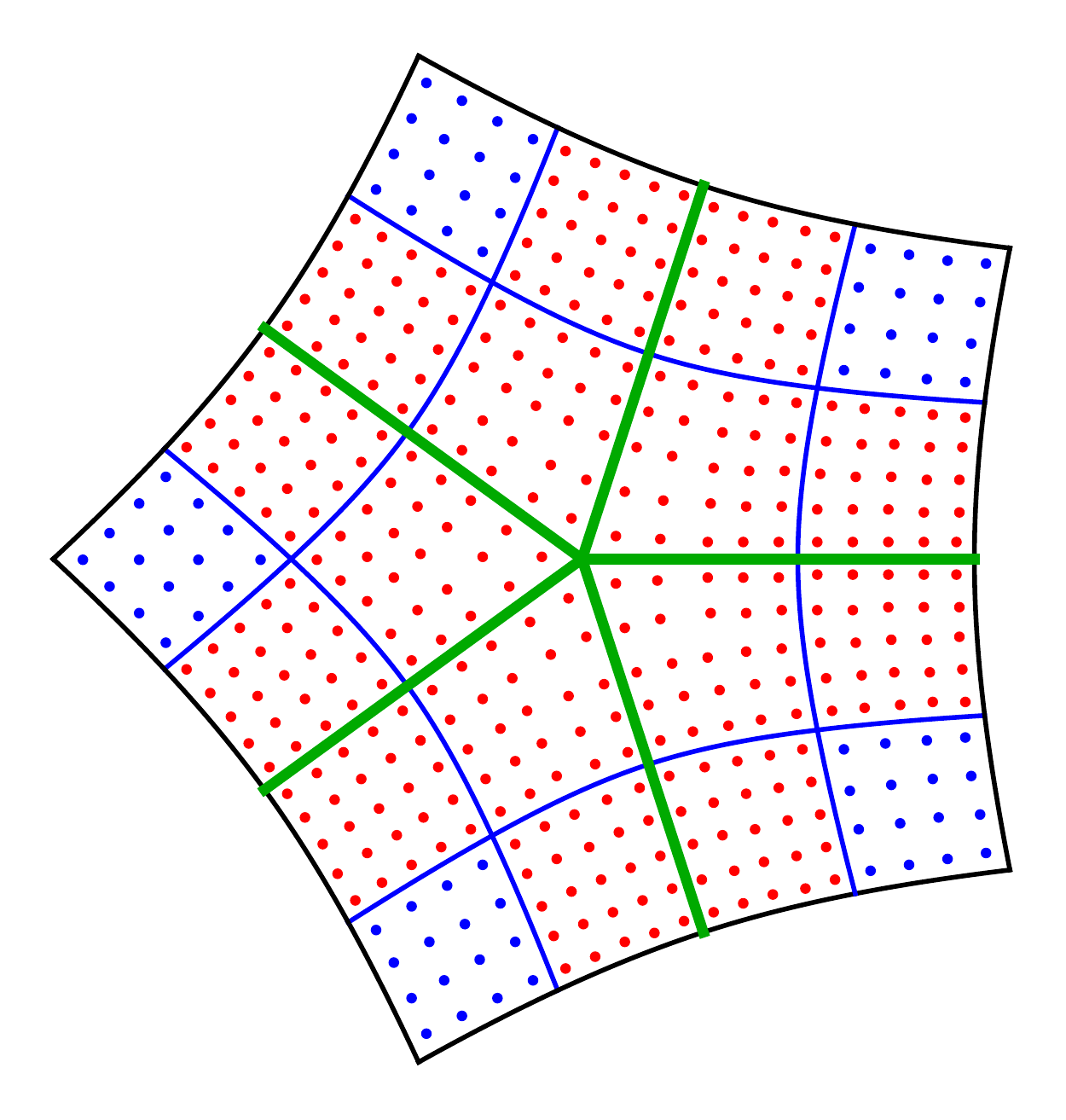}
\caption{A possible construction for an extraordinary vertex of valency $5$.}
\label{fig:ev5}
\end{figure}

The type of patches depicted in Figure \ref{subfig:2if} can be used to construct $C^1$ smooth isogeometric functions around extraordinary vertices. See Figure \ref{fig:ev5} for a 
possible construction.
We refer to \cite{Pe15-2,KaPe15A,KaPe18,NgKaPe15}, where such constructions were employed.

\begin{figure}[ht!]
\centering
\includegraphics[width=.2\textwidth]{4edges-ref0.pdf}\qquad $\xrightarrow{\mbox{refine}}$ \qquad
\includegraphics[width=.2\textwidth]{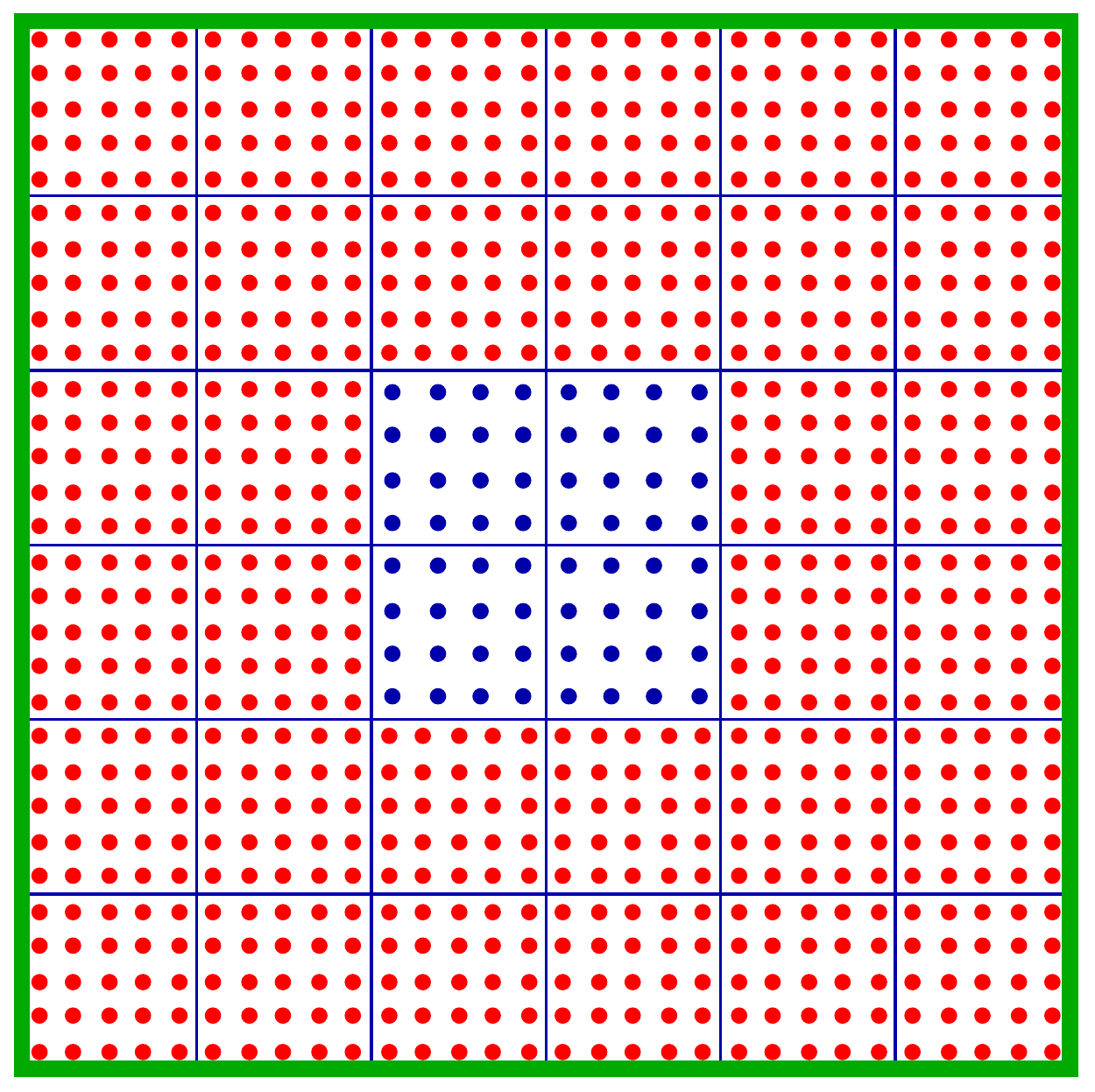}\qquad $\xrightarrow{\mbox{reduce}}$ \qquad
\includegraphics[width=.2\textwidth]{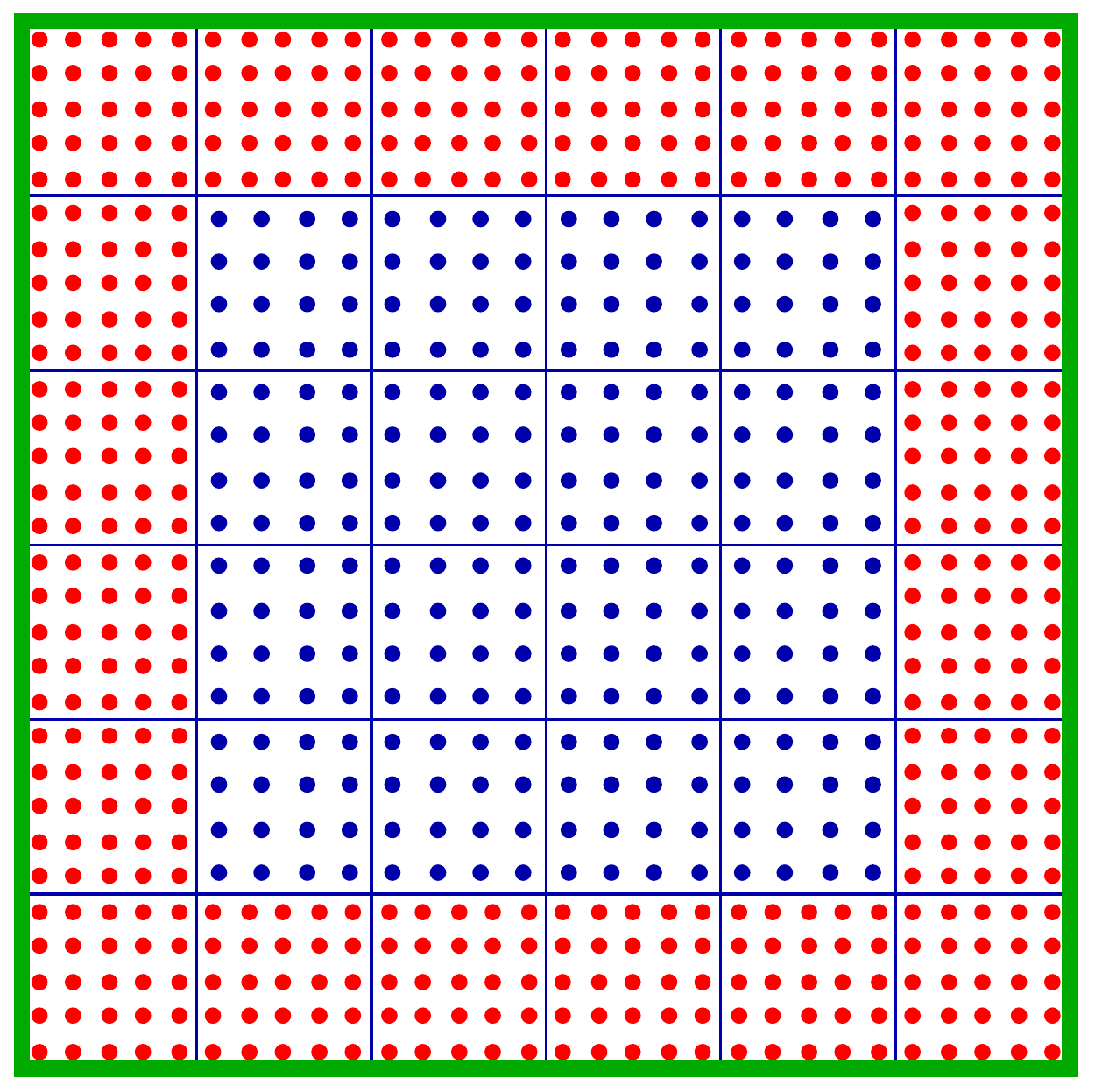}
\caption{A mixed degree patch (left), its refinement (middle) and a patch of reduced degree (right).}
\label{fig:multi-degree-patch-refinement}
\end{figure}

Another difficulty arises when refining the space. When performing a standard refinement step, the region where the degree is higher remains the same. Therefore the number of elements 
of higher degree scales with $O((\frac{1}{h})^2)$. This can be circumvented by locally reducing the degree again, which leads to the number of higher degree elements scaling as 
$O(\frac{1}{h})$. The process is sketched in Figure \ref{fig:multi-degree-patch-refinement}. Note that in this setting, the final (refined and reduced) space is not a superspace of the 
initial space. Hence, the spaces are not nested.

The constructions extend to higher degree as well as to more complex meshes. Many questions arise, that are worth to study in more detail; such as the definition of a basis forming a 
partition of unity, how to obtain nested spaces or how to efficiently construct domains and discretization spaces suitable for isogeometric analysis.

\section{Numerical examples} \label{sec:numerical_example}

\begin{figure}
\centering\footnotesize 
\begin{tabular}{cc}
\includegraphics[width=6.8cm,clip]{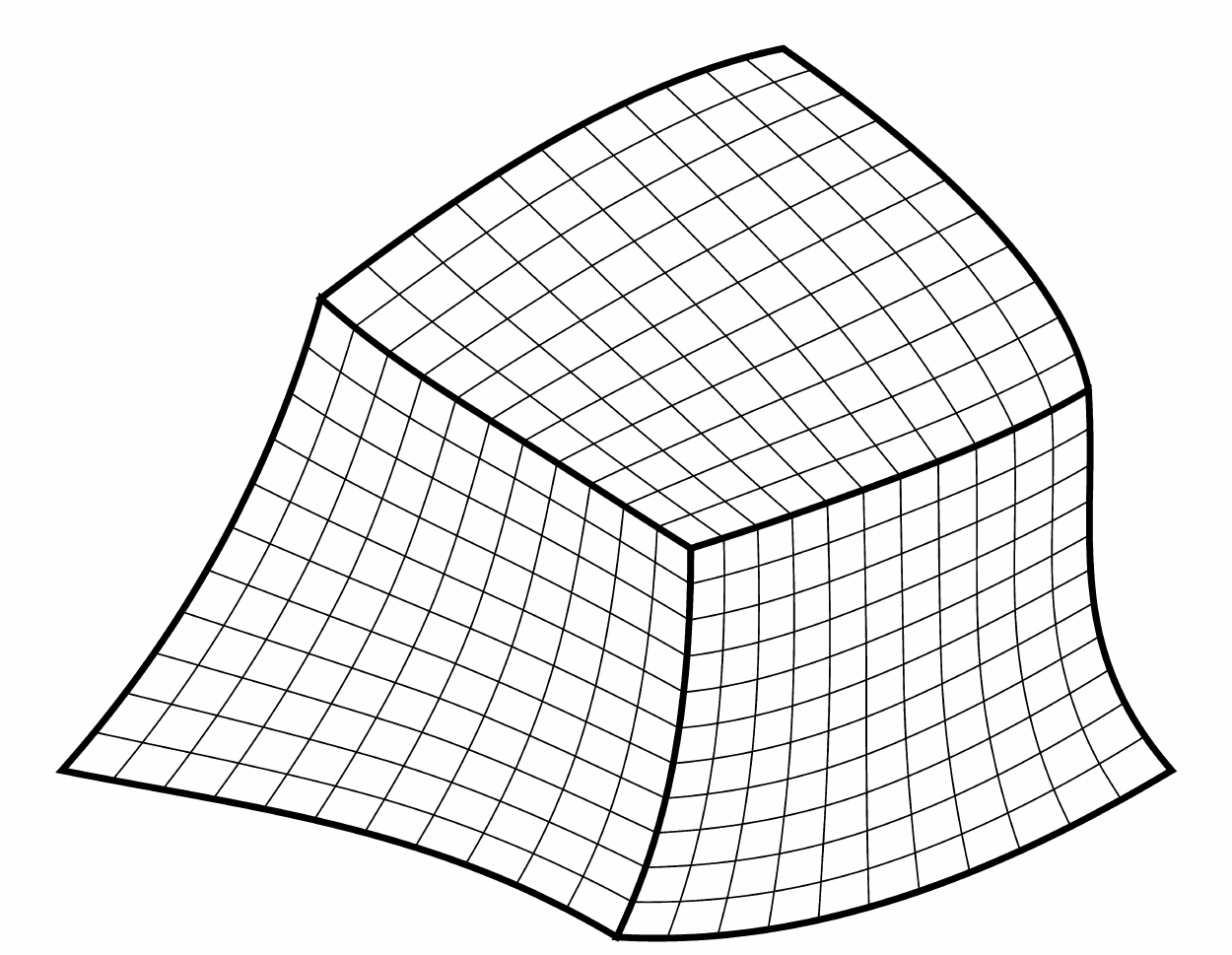} &
\includegraphics[width=6.2cm,clip]{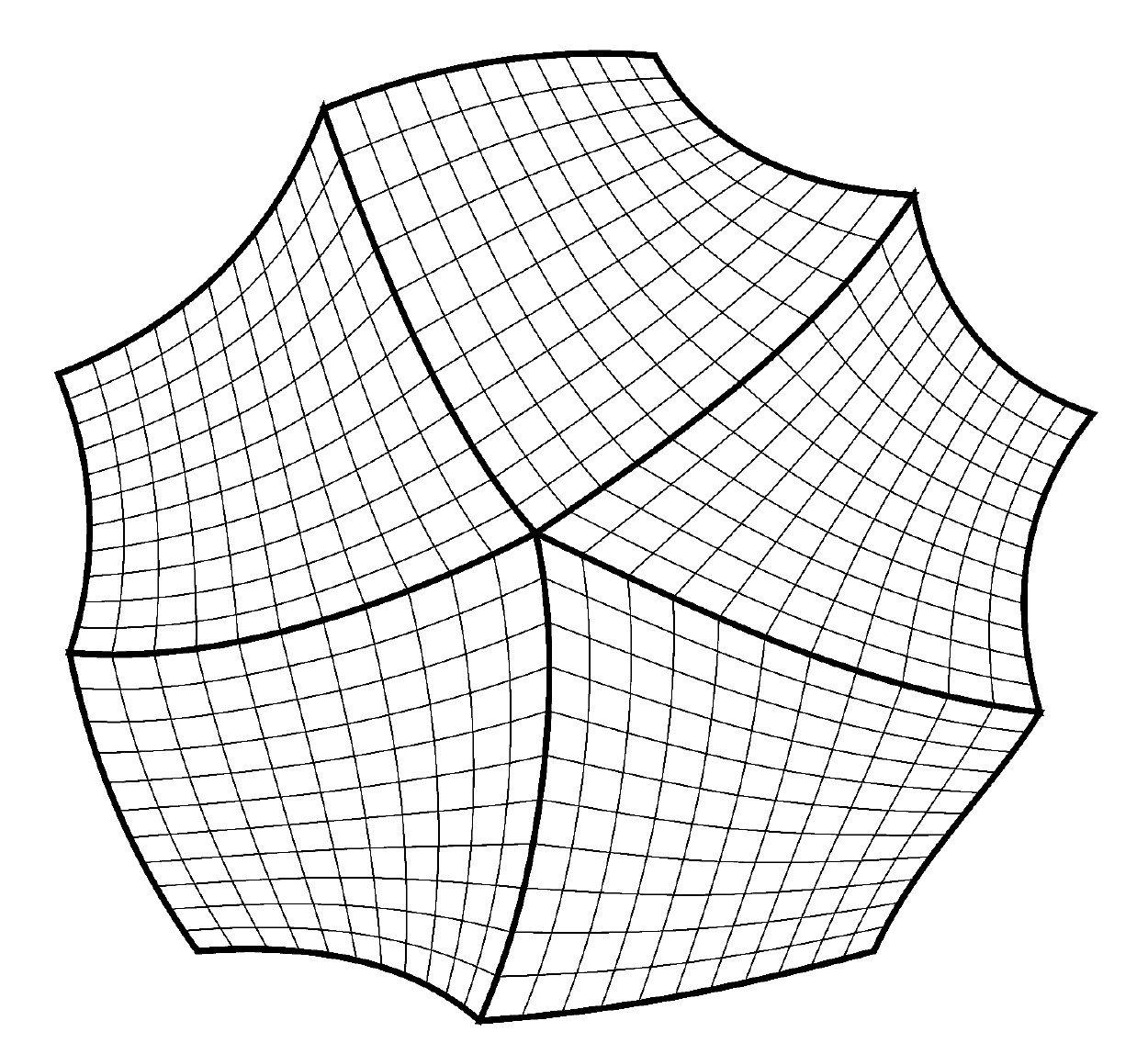} \\
\end{tabular}
\caption{Two AS-$G^1$ multi-patch geometries constructed my means of AS-$G^1$ fitting approach~\cite{KaSaTa17b}, cf. Section~\ref{subsec:approach_two}.}
\label{fig:AS_G1_geometries}
\end{figure}

We consider the two multi-patch domains~$\Omega$ shown in Fig.~\ref{fig:AS_G1_geometries}, which are described by AS-$G^1$ multi-patch geometries~$\f{F}$ consisting of 
parameterizations $\f{F}^{(i)} \in \left ( \mathcal{S}_{1/2}^{(3,3),(1,1)} \right )^2 $. The two AS-$G^1$ multi-patch geometries~$\f{F}$ have been constructed from initial multi-patch 
geometries~$\widetilde{\f{F}}$ composed of bicubic B\'{e}zier patches~$\widetilde{\f{F}}^{(i)}$ by using the AS-$G^1$ fitting approach~\cite{KaSaTa17b}, 
cf. Section~\ref{subsec:approach_two}. While the AS-$G^1$ three-patch geometry~(left) has been used in~\cite[Section~5]{KaSaTa17c}, too, the AS-$G^1$ five-patch 
parameterization~(right) is newly generated for this work.
For both multi-patch parameterizations~$\f{F}$ we generate a sequence of $C^1$ Argyris spaces~$\mathcal{A}_{h}$, $h=\frac{1}{4},\frac{1}{8},\frac{1}{16},\frac{1}{32}$, 
for $\f{p}=(p,p)=(3,3),(4,4)$ and $\f{r}=(r,r)=(1,1)$. 

\begin{figure}
\centering\footnotesize 
\begin{tabular}{ccc}
Exact solution & Relative errors for $p=3$ & Relative errors for $p=4$ \\  
\includegraphics[width=5.3cm,clip]{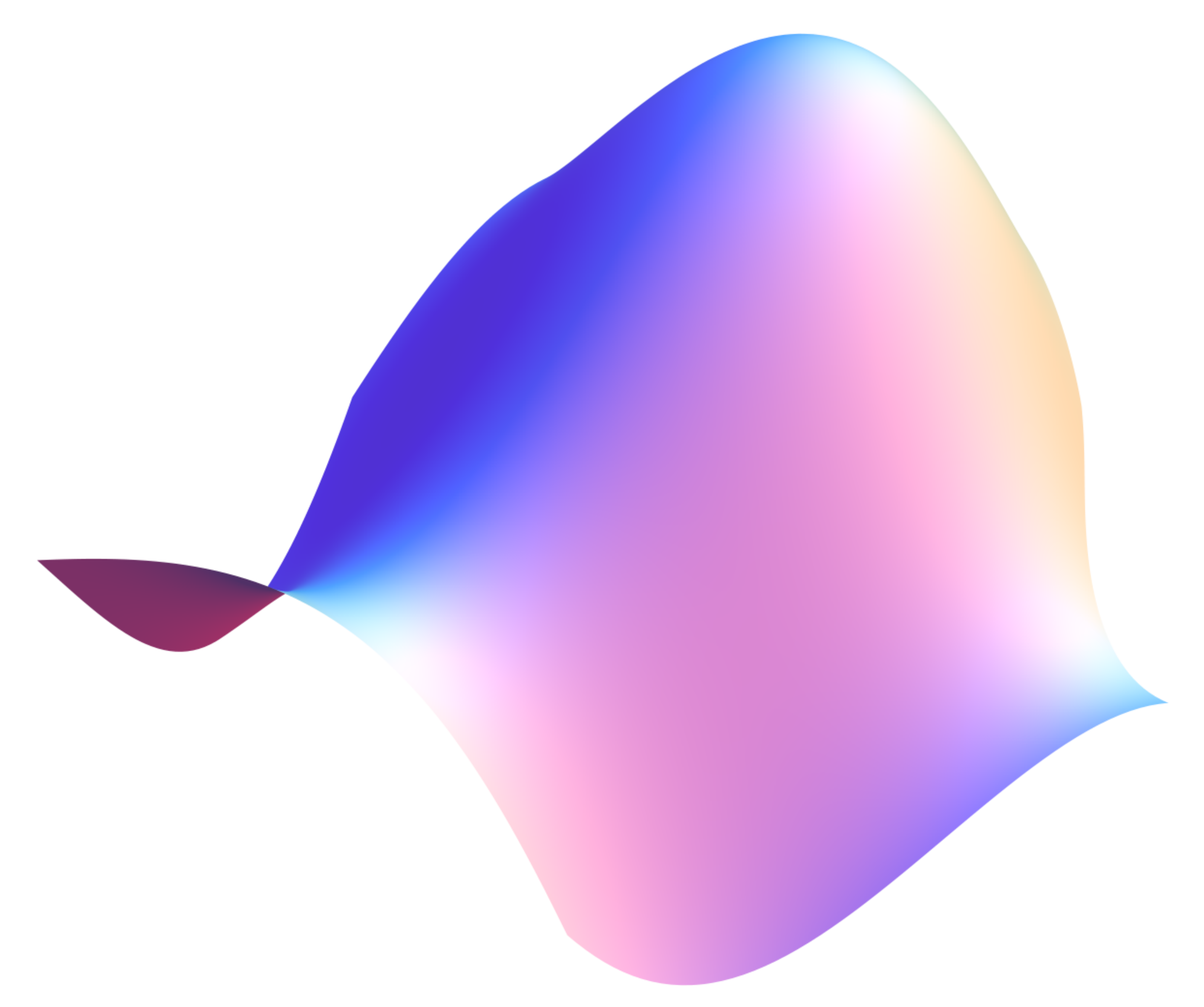} &
\includegraphics[width=5.0cm,clip]{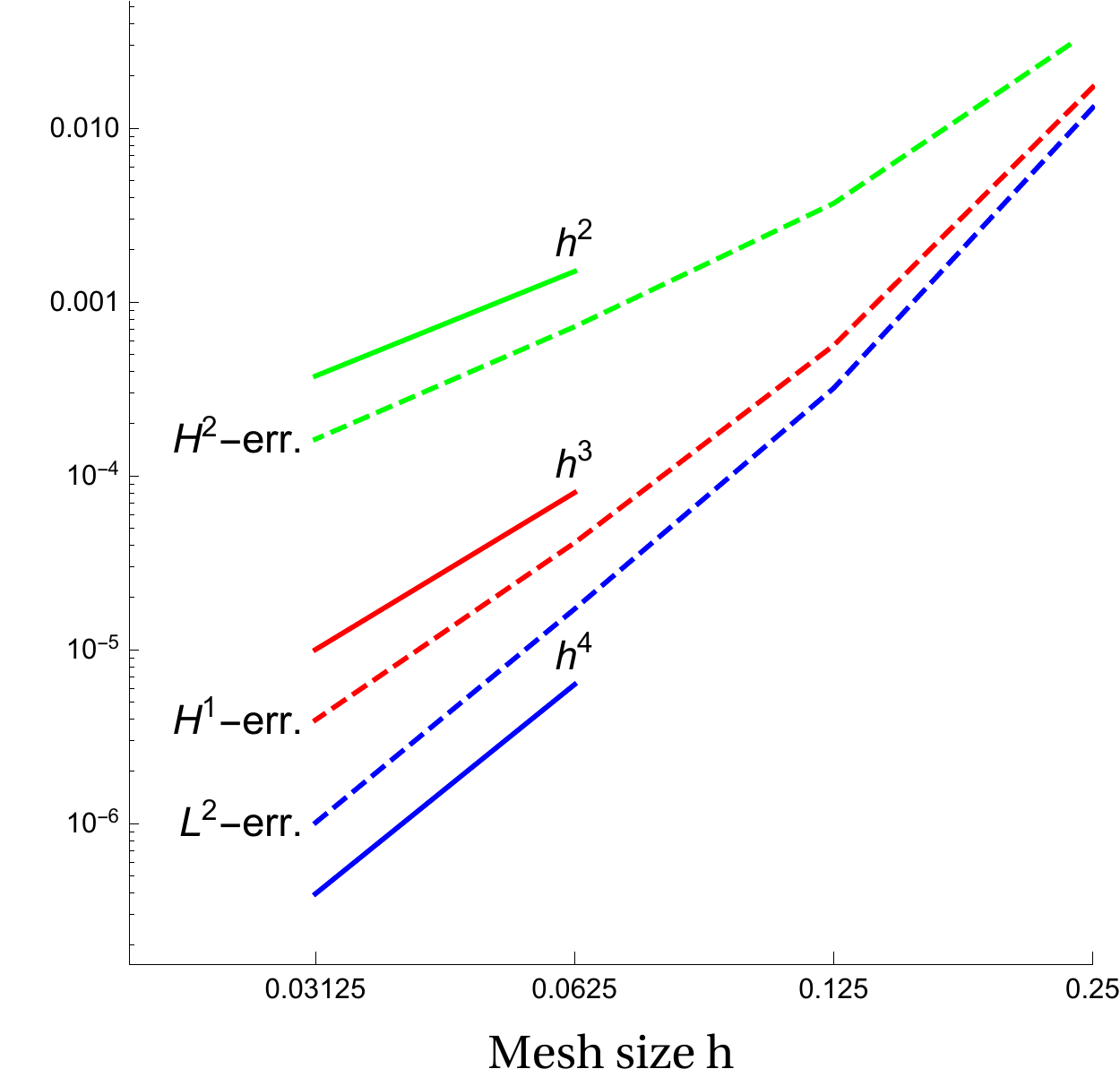} &
\includegraphics[width=5.0cm,clip]{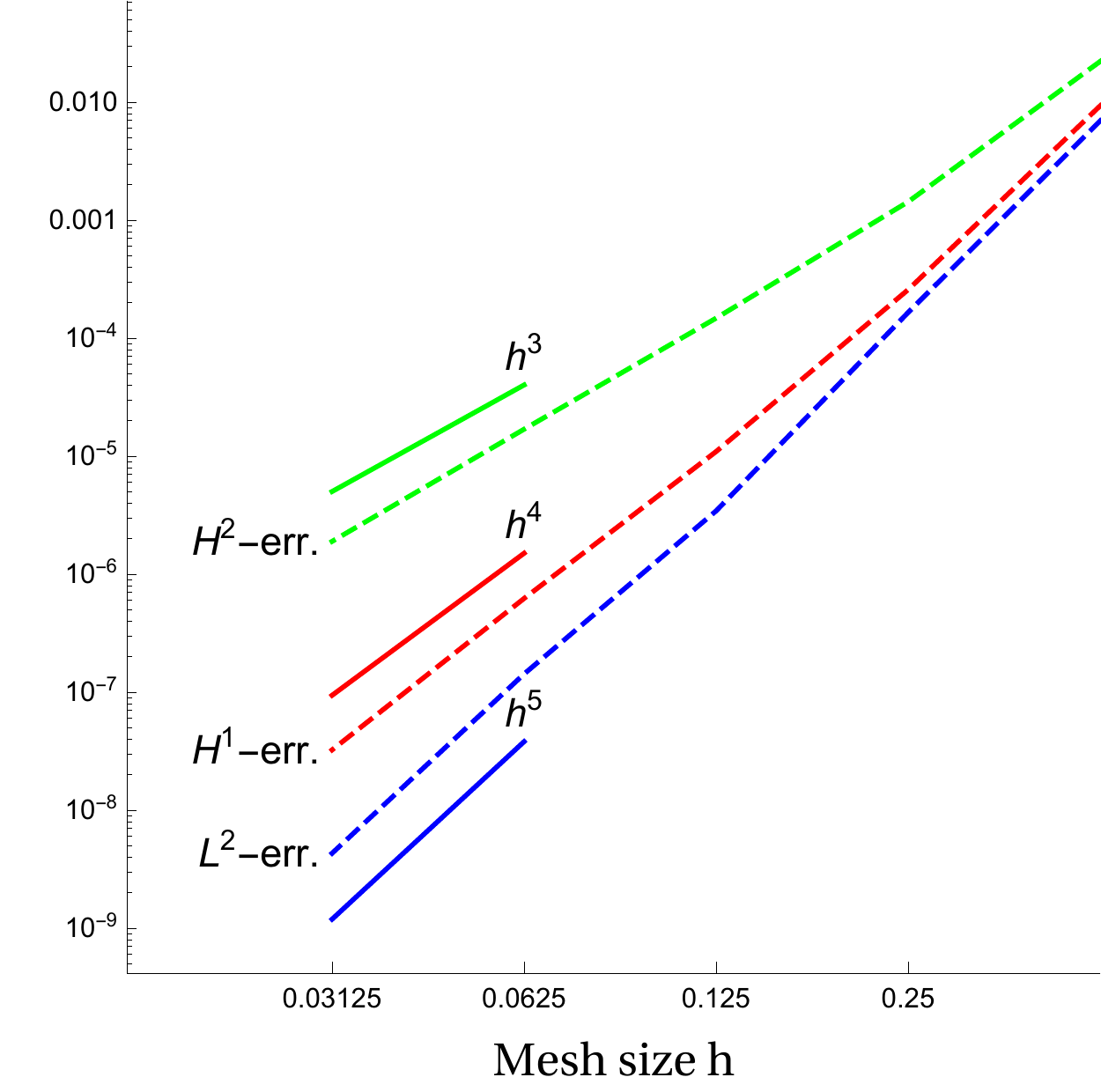} \\
\multicolumn{3}{c}{Example: AS-$G^1$ three-patch geometry} \\
\includegraphics[width=5.3cm,clip]{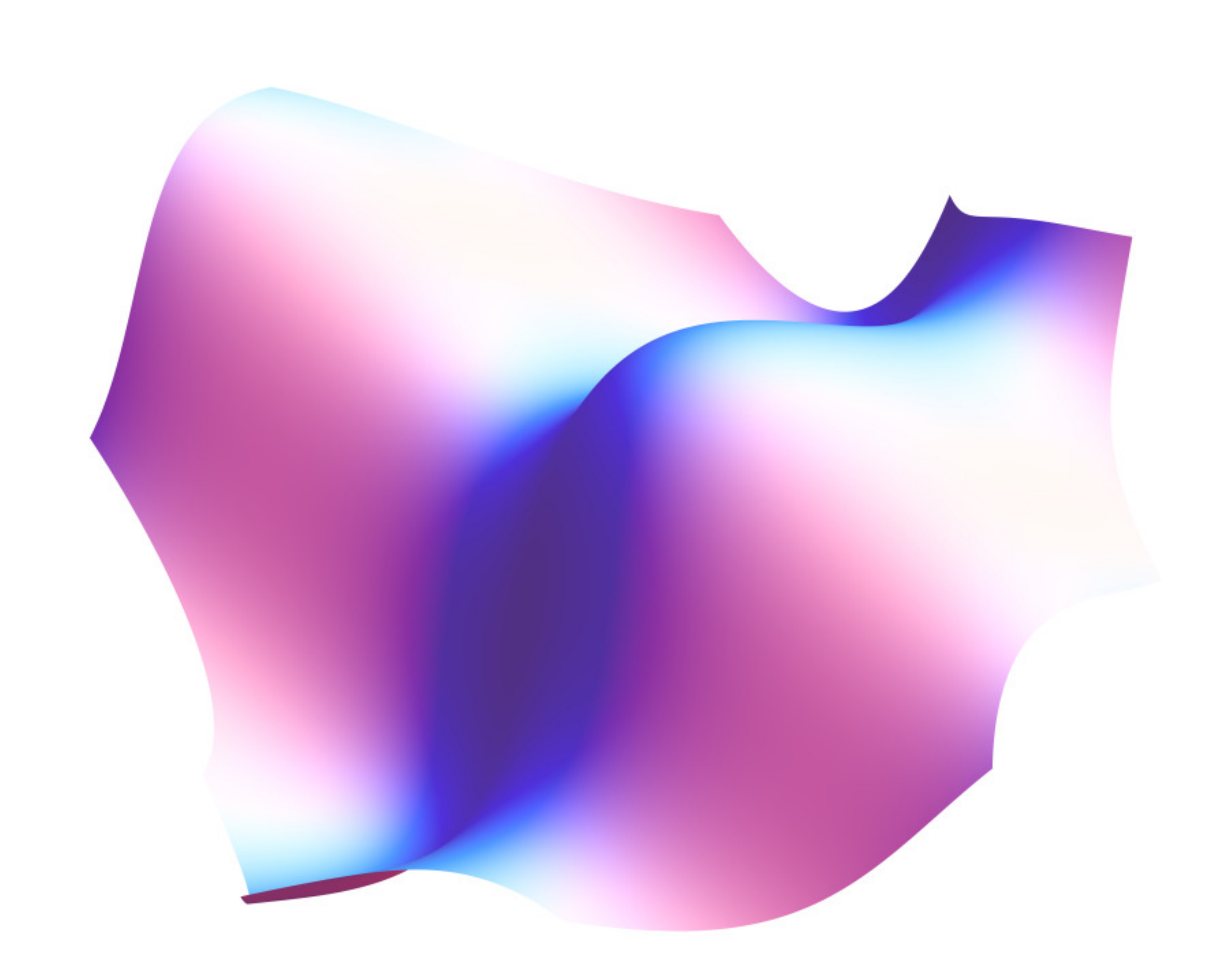} &
\includegraphics[width=5.0cm,clip]{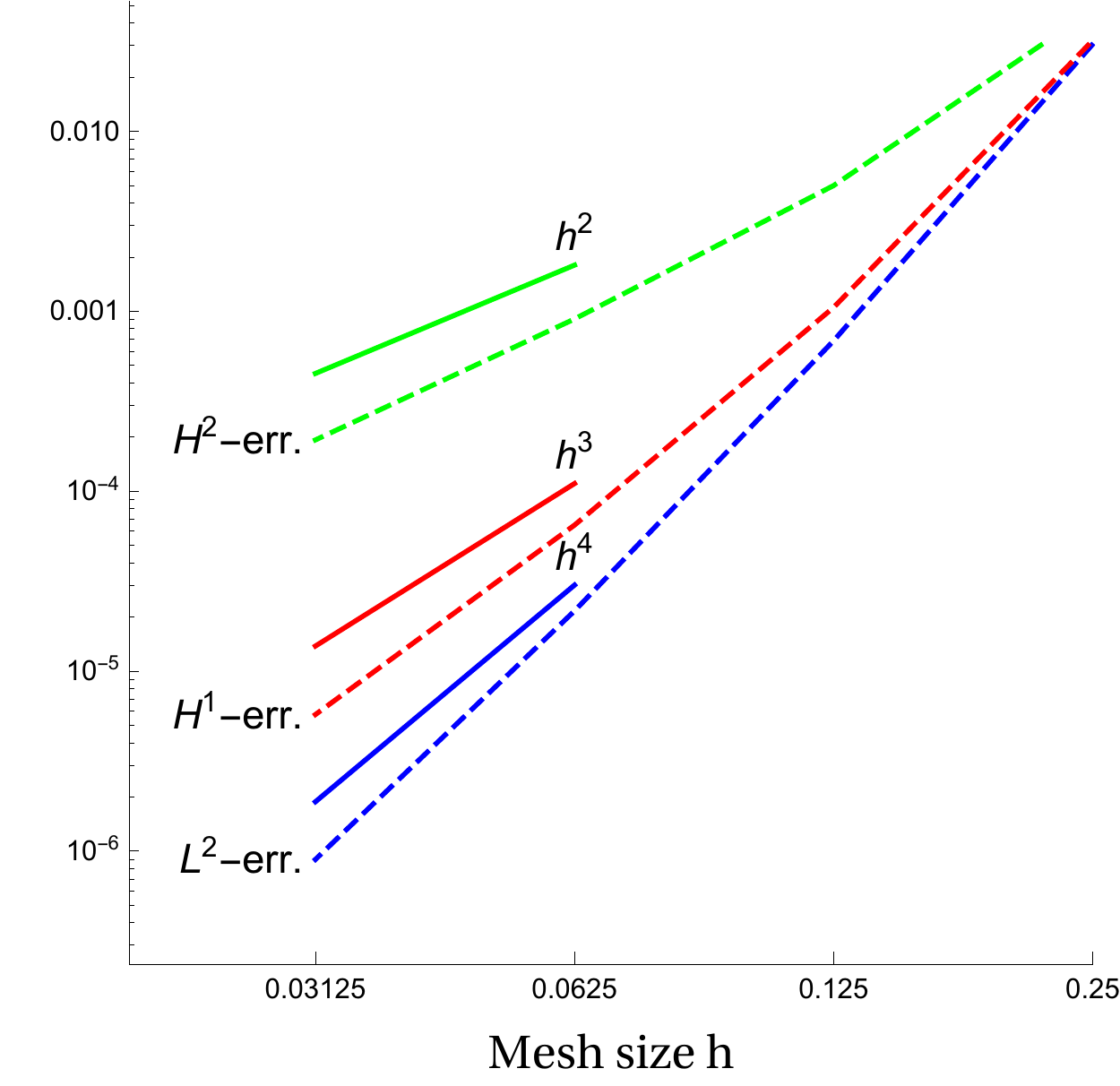} &
\includegraphics[width=5.0cm,clip]{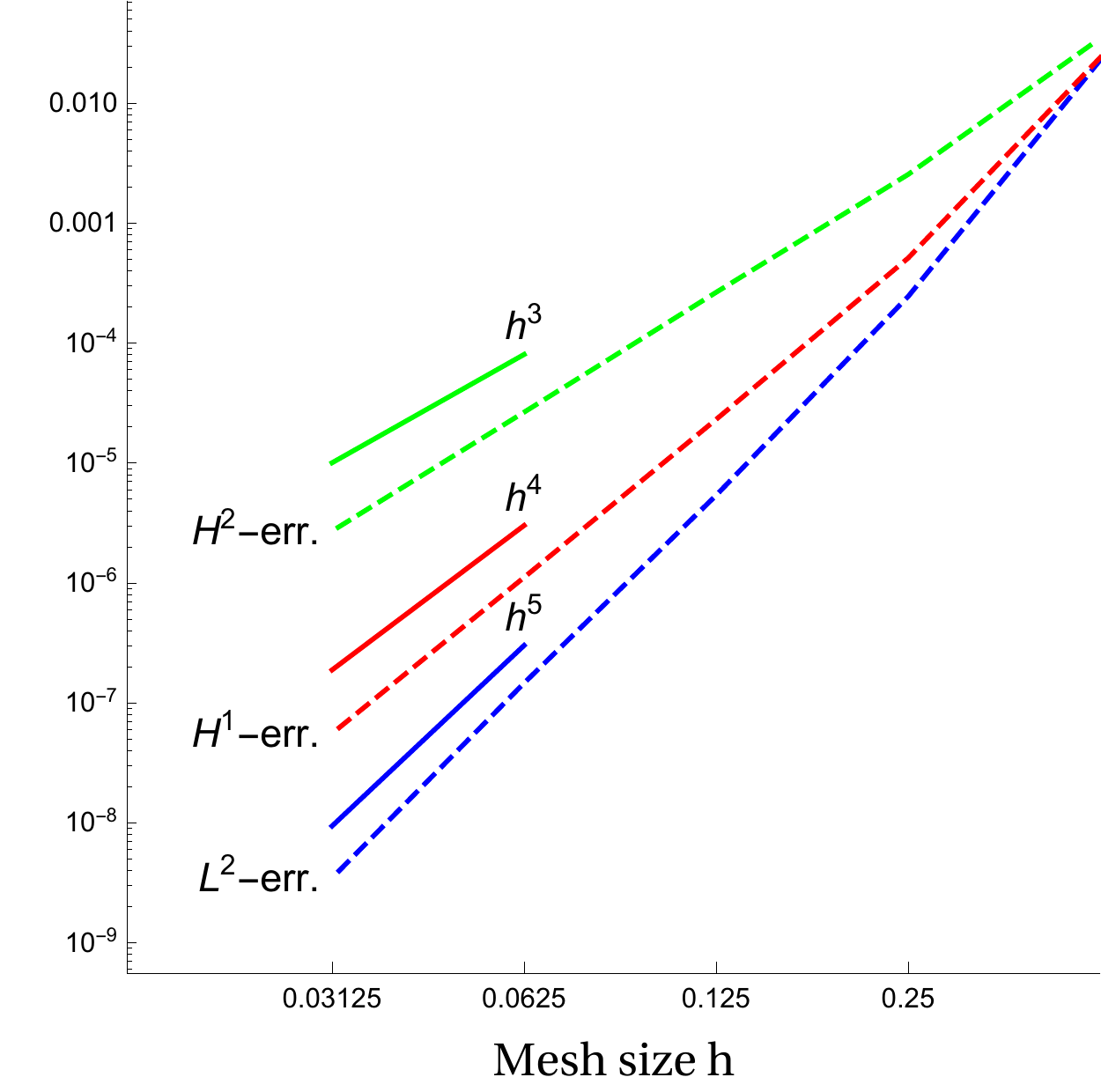} \\
\multicolumn{3}{c}{Example: AS-$G^1$ five-patch geometry}
\end{tabular}
\caption{Solving the biharmonic equation~\eqref{eq:problem_biharmonic} over the two AS-$G^1$ multi-patch geometries from Fig.~\ref{fig:AS_G1_geometries}: Exact solutions 
(first column) and the resulting relative $L^2$, $H^1$ and $H^2$ errors for $p=3$ (second column) and $p=4$ (third column).}
\label{fig:numerical_example}
\end{figure}

We employ the space family  $\mathcal{A}_{h}$ to solve the biharmonic equation
\begin{equation} \label{eq:problem_biharmonic}
\left\{
 \begin{array}{rll} 
 \triangle^{2} u(\f{x}) & = g(\f{x}) & \f{x} \in \Omega \\
  u(\f{x}) & = g_1(\f{x})  & \f{x} \in \partial \Omega \\
  \frac{\partial u}{\partial \f{n}}(\f{x})  & = g_2(\f{x}) &  \f{x} \in \partial \Omega
        \end{array} \right. 
 \end{equation}
by a standard  Galerkin discretization. The  
functions~$g$, $g_1$ and $g_2$ are selected to obtain the exact solution
\[
 u(\f{x}) = u (x_1,x_2)= -4 \cos (\frac{x_1}{2}) \sin (\frac{x_2}{2})
\]
on both domains.  In particular, the boundary Dirichlet data $g_1 $ and $g_2 $ are  $L^2$ projected and imposed strongly to the numerical solution $u_h$.  The resulting relative 
$L^2$, $H^1$ and $H^2$ errors with their estimated convergence rates are presented in Fig.~\ref{fig:numerical_example}~(second and third column), and indicate rates of optimal 
order of $\mathcal{O}(h^{p+1})$, $\mathcal{O}(h^{p})$ and $\mathcal{O}(h^{p-1})$, respectively.

\section{Conclusion}

In this paper we have listed and classified known methods to construct $C^1$-smooth isogeometric spaces over unstructured multi-patch domains. This is a research field that 
is attracting growing interest, at the confluence of geometric design and numerical analysis of partial differential equations.  We have discussed, with more details, the case of
multi-patch parametrizations that are regular and only $C^0$ at the patch interfaces, reviewing  in a coherent framework some of the recent  results that are  more closely related 
to our research activity.

\section*{Acknowledgments}
The research of G. Sangalli is partially supported by the European Research Council
through the FP7 Ideas Consolidator Grant \emph{HIGEOM} n.616563, and
by the Italian Ministry of Education, University and Research (MIUR)
through the  ``Dipartimenti di Eccellenza Program (2018-2022) - Dept. of Mathematics, University of Pavia''. 
The research of T. Takacs is partially supported by the Austrian Science Fund (FWF) and the government of Upper Austria through the project P~30926-NBL. 
This support is gratefully acknowledged. 
\bibliographystyle{plain}

\end{document}